\documentclass[11pt]{article}
\usepackage{graphicx}

\usepackage{style-paper}

\title{Quasipolynomial pyjamas}
\author{Sida Li\thanks{Trinity College, University of Cambridge, United Kingdom. Email: \href{mailto:sl2190@cam.ac.uk}{sl2190@cam.ac.uk}.}}
\date{}

\begin{document}
\maketitle

\begin{abstract}
    The ``pyjama stripe" $E(\eps)$ is the set of all complex numbers whose real part lies within $\eps$ of the nearest integer. The Pyjama Problem asks whether, for each $\eps > 0$, there exists some union of finitely many rotations of $E(\eps)$ about the origin that covers the plane. This was answered in the affirmative by Manners. Kravitz and Leng subsequently showed that $\exp(\exp(\exp(\eps^{-O(1)})))$ rotations suffice. In this paper, we prove an upper bound of
    \[\exp((\log \eps^{-1})^{50000})\]
    rotations. We also show that at least $(\pi/2 + 0.003 + o(1))\eps^{-1}$ rotations are necessary, improving upon the natural lower bound of $(\pi/2 + o(1))\eps^{-1}$.
\end{abstract}

\section{Introduction}

For $\eps > 0$, let $E(\eps)$ denote the ``pyjama stripe", given by
\[E(\eps) := \{z \in \mathbb{C} : \Re(z) \in (-\eps, \eps) \pmod{1}\}.\]
This is the central object of study in the Pyjama Problem of Iosevich, Kolountzakis and Matolcsi \cite{iosevich2006}. For $\theta$ a unit-norm complex number, let $E_\theta(\eps) := \{z \in \mathbb{C} : \Re(\theta z) \in (-\eps, \eps) \pmod{1}\}$ denote the pyjama stripe rotated by $\theta^{-1}$ about the origin. 

\begin{problem}\label{prob:pyjama}
    Let $\eps > 0$. Does there exist a finite set $\{\psi_1, \dots, \psi_N\}$ of unit-norm complex numbers such that
    \[\bigcup_{n = 1}^N E_{\psi_n}(\eps) = \mathbb{C}?\]
\end{problem}

For $\eps > 1/3$, one can easily verify that $E_1(\eps) \cup E_{e^{2\pi i/3}}(\eps) \cup E_{e^{4\pi i/3}}(\eps) = \mathbb{C}$. A natural next approach would be to consider roots of unity in general. However, Malikiosis, Matolcsi and Ruzsa \cite{malikiosis2016, malikiosis2013} demonstrated that for $\eps < 1/3$, rotations of $E(\eps)$ by finitely many roots of unity cannot cover $\mathbb{C}$.

In a breakthrough paper, Manners \cite{manners2013} resolved the Pyjama Problem in the affirmative for all $\eps > 0$. Motivated by Furstenberg's Theorem \cite{furstenberg1967} on $(\times 2, \times 3)$-invariant closed subsets of $\mathbb{R}/\mathbb{Z}$, Manners considers the semigroup generated by two multiplicatively independent rational rotations. In particular, he works with the unit-norm complex numbers given by
\[\theta_5 := \frac{1+2i}{1-2i}, \quad \theta_{13} := \frac{2+3i}{2-3i},\]
each labelled by the norm of the Gaussian prime in the numerator (or denominator). To resolve the Pyjama Problem, Manners adapts Furstenberg's Theorem to a specific solenoid, with actions by $\theta_5$ and $\theta_{13}$ taking the roles of multiplication by 2 and 3. 

Since Manners's proof appealed to compactness arguments, it provided no quantitative bound on the number of rotations required in terms of $\eps$. This was the subject of Problem 41 in Green's open problem list \cite{green2024}. In recent work \cite{kravitz2025}, Kravitz and Leng prove that
\[\exp(\exp(\exp(\eps^{-C})))\]
rotations suffice, for some constant $C > 0$. 

To obtain this bound, Kravitz and Leng perform a quantitative analysis of Manners's argument. A key input is their adaptation of a result of Bourgain, Lindenstrauss, Michel and Venkatesh \cite{bourgain2009}, which provides an effective version of a $(\times 2, \times 3)$ theorem of Rudolph \cite{rudolph1990} (a special case of Furstenberg's Conjecture, which generalises Furstenberg's Theorem). They also combined their careful analysis with new geometric ideas to quantify the remaining steps. 

While keeping the same blueprint as Kravitz and Leng, we consider actions corresponding to a growing number of primes, not just 5 and 13. We obtain the following bound. 

\begin{thm}\label{thm:pyjamas}
    For $\eps > 0$ sufficiently small, there exists a finite set $\Theta$ of unit-norm complex numbers such that
    \[|\Theta| \le \exp((\log \eps^{-1})^{50000})\]
    and $\bigcup_{\psi \in \Theta} E_\psi(\eps) = \mathbb{C}$. 
\end{thm}

We make no effort to optimise the exponent of $\log \eps^{-1}$, which currently relies on the best bounds for the zero-free strip for finite Hecke $L$-functions.

We also consider the lower bound for the Pyjama Problem. Kravitz and Leng \cite{kravitz2025} noted a trivial volume-counting lower bound of $\frac12 \eps^{-1}$ rotations. Moreover, by considering an arc length bound along the circle $S_{1-\eps}(0) := \{z \in \mathbb{C} : |z| = 1 - \eps\}$, we see that at least $(\frac{\pi}{2} + o(1))\eps^{-1}$ rotations are required\footnote{This was communicated to the author by Noah Kravitz and James Leng, and independently by Timothy Gowers.}. This provides a natural barrier to overcome. Here we provide a slight improvement.
    
\begin{thm}\label{thm:lower}
    Let $\eps > 0$ be sufficiently small. If $\Theta$ is a finite set of unit-norm complex numbers such that $\bigcup_{\psi \in \Theta} E_\psi(\eps) = \mathbb{C}$, then we must have
    \[|\Theta| \ge \left(\frac{\pi}{2} + \frac{1}{300}\right) \eps^{-1}.\]
\end{thm}

We achieve this by considering a small neighbourhood of $S_{2-\eps}(0)$ alongside $S_{1-\eps}(0)$. The arguments are purely geometric with no connection to techniques for the upper bound. It is likely that Theorem \ref{thm:lower} does not demonstrate the true behaviour of the lower bound, and so a natural next target would be a lower bound of the form $\gamma(\eps) \cdot \eps^{-1}$, where $\gamma(\eps) \to \infty$ as $\eps \to 0$. 

\section{Proof outline}

\subsection{Upper bound}

We start by explaining the setup of Manners's argument \cite{manners2013}. Then we describe the quantitative framework of Kravitz and Leng \cite{kravitz2025}. Finally, we move onto our improvements on their approach. 

\subsubsection{A dual problem.} Malikiosis, Matolcsi and Ruzsa \cite{malikiosis2013} noticed that the pyjama stripes have a natural dual interpretation. The continuous characters on $\mathbb{C}$ are
\[\chi_z : \mathbb{C} \to \mathbb{R}/\mathbb{Z}, \quad \chi_z: w \mapsto \Re(zw) \pmod{1},\]
for $z \in \mathbb{C}$. Thus, for $\theta$ a unit-norm complex number, we may equivalently write
\[E_\theta(\eps) = \{z \in \mathbb{C} : \chi_z(\theta) \in (-\eps, \eps)\}.\]
We can extend this definition to include all characters on $\mathbb{C}$, not just the continuous ones. These form a group $\widehat{\mathbb{C}}$, the Pontryagin dual of $\mathbb{C}$ (when equipped with the discrete topology). We may then define
\[E_\theta'(\eps) := \{\chi \in \widehat{\mathbb{C}} : \chi(\theta) \in (-\eps, \eps)\}.\]
Importantly, $\widehat{\mathbb{C}}$ is compact as the dual of a discrete space, with $E_\theta'(\eps)$ an open subset (under the dual topology given by pointwise evaluation). Using this, Malikiosis, Matolcsi and Ruzsa showed that Problem \ref{prob:pyjama} is equivalent to asking whether
\[\bigcup_{\psi \in S^1} E_\psi'(\eps) = \widehat{\mathbb{C}}\]
holds. However, as remarked by Manners, $\widehat{\mathbb{C}}$ is non-metrisable and generally difficult to work with. 

To get around this, Manners considers the semigroup $\Phi_\infty$ generated by the two rational rotations $\theta_5, \theta_{13} \in \mathbb{Q}(i) \cap S^1$, as described in the Introduction. Explicitly, we have
\[\Phi_\infty := \bigcup_{N \ge 1} \Phi_N, \quad \Phi_N := \{\theta_5^s\theta_{13}^t : 0 \le s, t \le N\}.\]
Manners noticed that by committing to this rotation set in advance, he only needs to evaluate a given character on elements of $\Phi_\infty$, instead of the whole of $S^1$. Thus, we can restrict $\mathbb{C}$ to a smaller subring before taking the dual formulation of the Pyjama Problem. The hope is that for a well-chosen subring, the dual space becomes more manageable than $\widehat{\mathbb{C}}$. 

The requirement is that this subring contains $\Phi_\infty$. This motivates the following choice:
\[\mathbb{A} := \mathbb{Z}[i] \left[\frac{1}{\overline{P_5}}, \frac{1}{\overline{P_{13}}}\right], \quad P_5 := 1+2i, \quad P_{13} := 2 + 3i.\]
As before, this has compact Pontryagin dual $\mathbb{X} := \widehat{\mathbb{A}}$, with an embedding of $\mathbb{C}$ as the continuous characters. In contrast to $\widehat{\mathbb{C}}$, $\mathbb{X}$ is now metrisable and can even be thought of as a quotient of $\mathbb{C} \times \mathbb{Q}_5 \times \mathbb{Q}_{13}$. 

The \emph{dual pyjama stripe} within $\mathbb{X}$ takes the form
\[E^*(\eps) := \{\chi \in \mathbb{X} : \chi(1) \in (-\eps, \eps)\}.\]
Using this, Manners considers an alternative perspective on the Pyjama Problem. Instead of rotating the dual pyjama stripe to cover $\mathbb{X}$, we equivalently fix $E^*(\eps)$ and consider ``rotating" points of $\mathbb{X}$ into $E^*(\eps)$. This motivates us to equip $\mathbb{X}$ with actions corresponding to rotations. For $\theta \in \Phi_\infty$ and $\chi \in \mathbb{X}$, let
\[(\theta \cdot \chi)(z) := \chi(\theta z).\]
This is a continuous action, so we can view $(\mathbb{X}, \Phi_\infty)$ as a topological dynamical system. 

\subsubsection{Rational rotations.} With this new perspective, we would be done if the orbit of every point in $\mathbb{X}$ intersects $E^*(\eps)$. It turns out this is too much to ask for, and this fact is manifested in the original geometric setting as follows.

Observe that $\Phi_\infty \subseteq \mathbb{Q}(i)$. Hence, the finite set $\Phi_N$ lies in the lattice $\frac{1}{D}\mathbb{Z}[i]$ for some non-zero $D \in \mathbb{Z}[i]$. Consider any $\theta \in \Phi_N$, so $D\theta \in \mathbb{Z}[i]$. Since $E(\eps)$ is invariant under translation by $\mathbb{Z}[i]$, the rotated pyjama stripe $E_\theta(\eps)$ is doubly-periodic with respect to $D\mathbb{Z}[i]$. 

One of the main results of Malikiosis, Matolcsi and Ruzsa \cite{malikiosis2013} is that these periodic rotations cannot cover all of $\mathbb{C}$ when $\eps < 1/3$. This illustrates the existence of ``rational obstructions": points with specific rational form which are guaranteed to be uncovered. To briefly demonstrate: suppose that $\Theta \subseteq \frac{1}{D}\mathbb{Z}[i]$, where $D$ is an odd integer ($D$ can always be chosen in this way). Each $\theta \in \Theta$ can be written as $\frac{r + si}{t}$ in reduced form, where $r, s, t$ are coprime integers. This forces $r-s, t$ to be odd, and a short computation shows that all points of the form $\theta \cdot \frac{1+i}{2}D$ have real part exactly 1/2. This means that $\frac{1+i}{2}D$ lies outside of the union of the rotated pyjama stripes. 

The key insight of Manners's proof is that an almost-converse holds: for sufficiently large $N$, $\mathbb{C} \setminus \bigcup_{\theta \in \Phi_N} E_\theta(\eps)$ is actually contained in the union of small balls around rational obstructions, which lie in a sparse lattice. This is known as a ``rationality lemma". Then he adds a few irrational rotations to cover the remaining set, adapting the ``irrational trick" trick of Malikiosis, Matolcsi and Ruzsa.

To establish the rationality lemma, Manners exploits the compactness of the dual space $\mathbb{X}$. He first shows that rational obstructions correspond to low-order torsion points in $\mathbb{X}$. Now recall that $\mathbb{X}$ has one complex and two $p$-adic components. Balls in the complex component, centred at these torsion points, also rotate amongst themselves under the action of $\Phi_\infty$. This is captured in the following infinitary rationality lemma, which forms the core of Manners's argument. 

\begin{lmm}\label{lmm:manners-rationality}
    If $Y \subseteq \mathbb{X}$ is a closed, $\Phi_\infty$-invariant set such that the image of the evaluation map $Y(1) = \{\chi(1) : \chi \in \mathbb{X}\}$ is not all of $\mathbb{R}/\mathbb{Z}$, then
    \[Y \subseteq \bigcup_{i = 1}^k \overline{B}_{1/2}^\mathbb{C}(q_i),\]
    for some finite set of torsion points $\{q_1, \dots, q_k\} \subseteq \mathbb{X}$. 
\end{lmm}

This lemma should be compared to the statement of Furstenberg's Theorem \cite{furstenberg1967}, which can be summarised as follows. If $S \subseteq \mathbb{R}/\mathbb{Z}$ is a closed, $(\times 2, \times 3)$-invariant set that is not the whole of $\mathbb{R}/\mathbb{Z}$, then $S$ is a finite set of rationals. However, there are noteworthy differences.

The most notable difference is that instead of simply assuming $Y \neq \mathbb{X}$, the rationality lemma uses the stronger assumption that $Y(1) \neq \mathbb{R}/\mathbb{Z}$. This encapsulates the ``notion of density" required for the Pyjama Problem as follows. By definition, the orbit of a non-torsion point is dense if it meets every non-empty open set. For a weaker notion, we call an orbit ``$\eps$-dense" if it meets every $\eps$-ball in $\mathbb{X}$ (with its prescribed metric). In the context of the Pyjama Problem, we can go even further. An orbit can be considered dense as long as it intersects $E^*(\eps)$, detected by evaluation at 1. This theme becomes crucial when we introduce a growing space.

The stronger assumption also allows Manners to constrain the radius of each complex ball, which could otherwise be arbitrarily large. We will discuss in more detail this later.

\subsubsection{The Kravitz-Leng framework.} \label{subsec:kravitz-leng} The main technical result of Kravitz and Leng is a quantitative rationality lemma of the following form. Take small parameters $\eps, \delta > 0$. Then for each point $p \in \mathbb{X}$ either: the orbit of $p$ under $\Phi_N$ intersects $E^*(\eps)$; or there is a torsion point $q$ of order at most $n$ such that $p$ is distance at most $\delta$ from $\overline{B}_{1/2}^\mathbb{C}(q)$. They demonstrate this with $n, N$ both triple-exponential in $\eps^{-1}$, and $N$ logarithmic in $\delta^{-1}$. 

In the proof sketch that follows, we introduce various parameters but do not keep track of the quantitative constraints until the end. We will then locate the sources of the three ``lost" exponentials.

The first step relies on the pigeonhole principle. They run an initial short orbit $\Phi_{N_0} \cdot p$ and pick out two points that are close together. In particular, there exist distinct $\phi_1, \phi_2 \in \Phi_{N_0}$ with
\[\dist((\phi_1 - \phi_2)p, 0) = \dist(\phi_1p, \phi_2p) \ll N_0^{-1/2}.\]
One can also bound the height of $r := \phi_1 - \phi_2 \in \mathbb{Q}(i)$, given by the maximum of the absolute values of the numerator and denominator (in reduced form). The proof proceeds by splitting into two cases: whether $rp$ is extremely close to $\overline{B}_1^\mathbb{C}(0)$, or merely somewhat close. 

The first case is labelled the ``major arcs". Since $r$ is a Gaussian rational, one expects to be able to ``divide through by $r$" to obtain that $p$ is still extremely close to a torsion point, now with order controlled by the height of $r$. However, the image of the complex ball $\overline{B}_1^\mathbb{C}(0)$ under this division can end up with radius $|r|^{-1}$, which may be significantly larger than the desired 1/2. 

Kravitz and Leng overcome this as follows. Suppose $p$ is indeed close to a complex circle of radius larger than 1/2. This circle itself already intersects $E^*(\eps)$ and maps to other such circles under the action by $\theta_5$. Motivated by this, they make geometric observations about $\mathbb{X}$ that force the orbit $\Phi_K \cdot p$ to intersect $E_1^*(\eps)$. This involves an application of Baker's Theorem from Diophantine approximation \cite{bakerI1967, bakerII1967, bakerIII1967}. 

The other case is labelled ``minor arcs". This is when $rp$ is not extremely close to $\overline{B}_1^\mathbb{C}(0)$. By a quantitative analysis of Manners's proof of Lemma \ref{lmm:manners-rationality}, they demonstrate that a further long orbit of the difference set
\[\Phi_{N_1} \cdot (\Phi_{N_0} \cdot p - \Phi_{N_0} \cdot p) = \Phi_{N_0 + N_1} \cdot p - \Phi_{N_0 + N_1} \cdot p\]
is $\eta$-dense in $\mathbb{X}$. We call this the \emph{first mixing stage}. 

For the \emph{second mixing stage}, we wish to convert from an $\eta$-dense difference set to an $\eps$-dense orbit. Since evaluation at 1 is 1-Lipschitz, this is enough for the orbit to intersect $E^*(\eps)$. 

Kravitz and Leng construct a fine partition $\mathcal{P}$ of $\mathbb{X}$ such that the difference set of $\Phi_{N_0 + N_1} \cdot p$ intersects every part. From their construction, they deduce that the orbit itself intersects $\gg |\mathcal{P}|^{1/2}$ parts of $\mathcal{P}$. This condition becomes the input for the entropy machinery of Bourgain, Lindenstrauss, Michel and Venkatesh \cite{bourgain2009}, which Kravitz and Leng adapt to $\mathbb{X}$. The output is that a slightly longer orbit $\Phi_{N_0 + N_1 + N_2} \cdot p$ is $\eps$-dense in $\mathbb{X}$, as desired. 

By tracing through the quantitative dependencies, they obtain the following constraints:
\[\eta := 5^{-\eps^{-C_1}}, \quad N_2 := \eta^{-0.9}, \quad N_0 = \exp(\exp(\eps^{-C_2})) \gg 13^{-2C_1\eta^{-5}},\]
\[K := \exp(200C_3N_0), \quad n := \exp(200N_0), \quad N_1 := \max\{\log(\delta^{-1}) + \log n, 10K\} + C_1\eta^{-5},\]
where $C_1, C_2, C_3 > 0$ are constants.

We first lose an exponential in the choice of $\eta$, to make the entropy machinery work in the second mixing stage. We then lose another within the first mixing stage to define $N_0$. The last exponential is lost due to the pigeonhole step. In particular, the maximum height of $r$ is exponential in the size of the initial orbit, which is reflected in both $K$ and $n$. 

\subsubsection{Our improvements.} \label{subsec:improvements} The source of the first two exponentials can be summarised as follows. The set $\Phi_N$ has cardinality roughly $N^2$, but maximum 5-adic and 13-adic absolute value both exponential in $N$. This results in quantitative losses because the metric on $\mathbb{X}$ captures these $p$-adic absolute values.

Motivated by this, it becomes natural to consider spreading out orbits across many primes, i.e. actions by $\theta_p$ for rational primes $p$ other than just 5 and 13. Since we wish to place stronger bounds on $p$-adic absolute values, we restrict the exponents of these primes, and so we take the number of primes to grow with $\eps^{-1}$. 

Thus, given suitable primes $p_1, \dots, p_k$, the desired number theoretic input to both mixing stages is as follows. We wish to cover the multiplicative group modulo $P_p^n$ using short products of $\theta_{p_1}, \dots, \theta_{p_k}$, where $P_p$ is a Gaussian prime of norm $p$ and $n$ is a natural number. This gives us a density result on a discretised scale. As such, we call $p$ the \emph{scaling prime} and $p_1, \dots, p_k$ the \emph{mixing primes}. 

This number theoretic result draws inspiration from a recent work of Klurman, Shparlinski and Teräväinen \cite{klurman2026}, related to a paper of Pilatte \cite{pilatte2026}. The main theorem of \cite{klurman2026} is the following. Given a parameter $\eps > 0$ and a threshold $Q \ge 3$, one can find a set of at most $(\log Q)^{1+\eps}$ primes of polylogarithmic size whose subset products generate $(\mathbb{Z}/q\mathbb{Z})^\times$, for almost all moduli $q \le Q$. The noteworthy aspects of this statement are the complete coverage of coprime residues, and the number of subset products is only quasipolynomial in $Q$. 

To prove such a result, we adapt the techniques of Klurman, Shparlinski and Teräväinen to finite Hecke $L$-functions, obtaining power-saving bounds for a character sum over Gaussian primes. We refer to Section \ref{subsec:hecke} for an introduction to the relevant theory.\\

A growing number of primes introduces a new obstruction, which relates to our above discussion on notions of density. The appropriate dual space $X$ is now a quotient of $\mathbb{C} \times \mathbb{Q}_5 \times \mathbb{Q}_{13} \times \cdots \times \mathbb{Q}_{p_k}$. Thus, one can imagine that the $\eps$-ball has volume roughly $\eps^{k+2}$, so every $\eps$-net of $X$ has size at least $\eps^{-k}$. If $k$ grows at least as fast as $\eps^{-1}$, an orbit of quasipolynomial size cannot possibly become $\eps$-dense within $X$.

To bypass this, we adapt the minor arcs argument of Kravitz and Leng with refined notions of density at each mixing stage. We first weaken the final $\eps$-density statement to pyjama density (intersecting $E^*(\eps)$). The following critical observation allows us to do this. In complex translates of the standard fundamental domain of $X$, given by
\[([-1/2, 1/2) \times [-1/2, 1/2)i) \times \mathbb{Z}_5 \times \cdots \times \mathbb{Z}_{p_k},\]
membership in the dual pyjama stripe is determined by only the real coordinate. 

Our entropic proof in the second mixing stage can then be catered directly to $E^*(\eps)$, which more or less lets us ``ignore" the $p$-adic components. This heavily simplifies the arguments of Kravitz and Leng \cite{kravitz2025}, practically to the original formulation in $\mathbb{R}/\mathbb{Z}$ by Bourgain, Lindenstrauss, Michel and Venkatesh \cite{bourgain2009}. 

Working backwards, we also choose the output of our first mixing stage to coincide exactly with the level of density necessary for the second mixing stage. In particular, we relax $\eta$-density, and now we only require the difference set to intersect every part of a given partition of $X$. 

To shave the last exponential, we observe that the rotations used for the pigeonhole step can be different from the rotations used in the rest of the rationality lemma. This allows us to pick a more height-friendly set of rational rotations for this initial orbit. 

\subsection{Organisation of the paper}

We split the paper between our upper bound and lower bound. In particular, we prove Theorem \ref{thm:pyjamas} in Section \ref{upper}, and then we prove Theorem \ref{thm:lower} in Section \ref{lower}.

In Section \ref{sec:prelim}, we build the background required for our proof. This includes properties of our space $X$ in Section \ref{subsec:X}, constructing our fundamental domains in Section \ref{subsec:fundamental}, an overview of the group $(\theta_5^{-n}A)/A$ in Section \ref{subsec:theta5}, and an introduction to the relevant theory of finite Hecke characters in Section \ref{subsec:hecke}. 

Then we proceed with our modification of the Kravitz-Leng framework. We cover the initial pigeonhole step, major arcs, minor arcs and the conclusion in Sections \ref{sec:pigeon}, \ref{sec:major}, \ref{sec:minor}, \ref{sec:conversion} respectively. 

Section \ref{lower} covers lower bounds. We start with a detailed treatment of the circle lower bound in Section \ref{sec:first-circle}, as input to our proof of the improved lower bound in Section \ref{sec:second-circle}. 

\section{Upper bound} \label{upper}

\subsection{Preliminaries}\label{sec:prelim}

In this section, we develop the theory used to motivate the proof. 

\subsubsection{The growing space $X$.}\label{subsec:X} We work with $k$ Gaussian primes, where $k$ is a positive integer parameter. We adopt the following conventions. 

Let $P_p$ denote the Gaussian prime with $N(P_p) = p$, where $p \equiv 1 \pmod{4}$ is a rational prime, satisfying $P_p = a + bi$ where $b > a > 0$. For instance, continuing the examples in the Introduction,
\[P_5 := 1 + 2i, \quad P_{13} := 2 + 3i, \quad P_{17} := 1+4i, \quad \dots\]
We consider products of actions of the form $\theta_p := P_p/\overline{P_p}$, noting that $|\theta_p| = 1$ and $\theta_p \in \mathbb{Q}(i)$. In particular, they all lie in the ring
\[A := \mathbb{Z}[i] \left[\frac{1}{\overline{P_5}}, \frac{1}{\overline{P_{13}}}, \dots, \frac{1}{\overline{P_{p_k}}}\right],\]
where $p_1, p_2, \dots, p_k$ denote the smallest $k$ primes such that $p_j \equiv 1 \pmod{4}$. 

Let $X := \widehat{A}$ denote the Pontryagin dual of $A$ (when equipped with the discrete topology). We seek a more explicit description of this dual space. Here we draw analogues to Chapter 3.2 of \cite{kravitz2025}, which covers the case of $A$ replaced by
\[\mathbb{Z}[i]\left[\frac{1}{\overline{P_5}}, \frac{1}{\overline{P_{13}}}\right],\]
and to Proposition 6.2 of \cite{manners2013}, which covers the same case after working in a model of $\mathbb{Z}[1/2]$. 

First, we define the dual pyjama stripe. For $\theta \in A$, let
\[E_\theta^*(\eps) := \{\chi \in X : \chi(\theta) \in (-\eps, \eps)\}, \quad E^*(\eps) := E_1^*(\eps).\]

For each prime $p_j$, we can complete $\mathbb{Q}(i)$ with respect to the $\overline{P_{p_j}}$-adic metric to obtain $\mathbb{Q}_{p_j}$. For instance, let $\iota_5: \mathbb{Q}(i) \to \mathbb{Q}_5$ denote the embedding which sends $i$ to the square root of $-1$ that is $-2 \pmod{5}$. Thus, we obtain that
\[|\iota_5(P_5)|_5 = |P_5|_{\overline{P_5}} = 1, \quad |\iota_5(\overline{P_5})|_5 = |\overline{P_5}|_{\overline{P_5}} = 5^{-1}.\]
We can identically repeat for $p_1, \dots, p_k$. We also embed $\mathbb{Q}(i)$ into $\mathbb{C}$ the usual way and drop the $\iota_\mathbb{C}$ for the sake of convenience. 

Recall that the integer and fractional parts of an element $z \in \sum_{l = u}^\infty a_l 5^l \in \mathbb{Q}_5$ are given by
\[\floor{z}_5 := \sum_{l = 0}^\infty a_l5^l, \quad \{z\}_5 := z - \floor{z}_5 \in [0,1),\]
respectively. Thus, we can define, for $x \in A$,
\[j_5:\mathbb{Q}_5 \to X, \quad j_5(a) : x \mapsto \{a\iota_5(x)\}_5.\]
We can similarly repeat for $p_1, \dots, p_k$, and we also introduce
\[j_\mathbb{C} : \mathbb{C} \to X, \quad j_\mathbb{C}(z) : x \mapsto \Re(zx) \pmod{1}.\]
Combining these maps altogether, we obtain
\[j:\mathbb{C} \times \mathbb{Q}_5 \times \cdots \times \mathbb{Q}_{p_k} \to X\]
\[j(z, a_5, \dots, a_{p_k}) := -j_\mathbb{C}(z) + j_5(a_5) + \cdots + j_{p_j}(a_{p_k}).\]
Chapter 3 of Chothi, Everest and Ward \cite{chothi1997} demonstrates that the kernel of $j$ is the diagonal embedding of $A$. This map is given by
\[\iota^\Delta : \mathbb{Q}(i) \to \mathbb{C} \times \mathbb{Q}_5 \times \cdots \times \mathbb{Q}_{p_k},\]
\[\iota^\Delta(z) := (z, \iota_5(z), \dots, \iota_{p_k}(z)).\]
Furthermore, we have that $j$ surjects onto $X$. Thus, $j$ induces an isomorphism of topological groups
\[(\mathbb{C} \times \mathbb{Q}_5 \times \cdots \times \mathbb{Q}_{p_k})/\iota^\Delta(A) \to X.\]
As the Pontryagin dual of a discrete space, $X$ is a compact space. Furthermore, it is endowed with the translation-invariant metric
\[\dist(x, y) := \inf\left\{|z|_\mathbb{C} + \sum_{j = 1}^k |a_{p_j}|_{p_j} : j(z, a_5, \dots, a_{p_k}) = x - y\right\}.\]
We also consider balls just in the complex component. For $R > 0$ and $p \in X$, let
\[\overline{B}_R^\mathbb{C}(p) := p + j_\mathbb{C}(\overline{B}_R(0)), \quad S_R^\mathbb{C}(p) := p + j_\mathbb{C}(S_R(0))\]
denote the complex ball and circle of radius $R$ around $p$, respectively. 

It follows from the above isomorphism that the torsion points of $X$ are $j(\iota^\Delta(\mathbb{Q}(i)))$, with the order of $q = j(\iota^\Delta(z))$ being the smallest positive integer $n$ such that
\[nq = 0 \in X \iff nz \in A.\]
Generalising Proposition 6.4 of \cite{manners2013}, these are also exactly the points with finite orbit under actions by all the $\theta_{p_j}$.

\subsubsection{Fundamental domains.} \label{subsec:fundamental} Our main focus is the standard fundamental domain of $X$,
\[F_0 := ([-1/2,1/2) \times [-1/2,1/2)i))\times \mathbb{Z}_5 \times \cdots \times \mathbb{Z}_{p_k}.\]
To demonstrate its utility, if $\chi = j(z, a_5, \dots, a_{p_k})$ is a representative of a character within $F_0$, then for $w \in A$,
\[\chi(w) = -\Re(zw) + \sum_{j = 1}^k \{a_{p_j}\iota_{p_j}(w)\}_{p_j} \pmod{1}.\]
In particular, since $a_{p_j}, \iota_{p_j}(1) \in \mathbb{Z}_{p_j}$, we have that $\chi(1) = -\Re(z) \pmod{1}$. Thus, our dual pyjama stripe has the following ``prism-like" shape within $F_0$,
\[E^*(\eps) = (E(\eps) \cap ([-1/2,1/2) \times [-1/2,1/2)i)) \times \mathbb{Z}_5 \times \cdots \times \mathbb{Z}_{p_k} \subseteq F_0.\]
We remark that the choice of $[-1/2, 1/2)$ instead of $[0,1)$ is simply notational, in order to make our choice of smooth bump function in Section \ref{subsec:entropy} more natural.

As introduced in Section \ref{subsec:improvements}, we consider one scaling prime and many mixing primes. For simplicity, we take 5 to be our scaling prime in most statements. The generalisations to other primes $p_j$ are natural and are encapsulated in the leeway in our bounds. 

Let $n$ be a non-negative integer with $N$ satisfying $N := 5^n$. In preparation for our partition on $X$, we look for fundamental domains of $X/(\theta_5^{-n}A)$. 

Let $h, \alpha$ be integers with $2 \le h \le k$. The intuition is that the primes $p_1, \dots, p_h$ are those used in both mixing stages (see Section \ref{subsec:manners} and Section \ref{subsec:entropy}), whereas the primes $p_{h+1}, \dots, p_k$ are just used in the first. 

With these in mind, consider the following subset of $X$,
\[F(n, h,k,\alpha) := \left(P_5^{-n}\overline{P_5}^\alpha \prod_{j = 2}^h \overline{P_{p_j}}\right) ([-1/2,1/2) \times [-1/2,1/2)i)) \times 5^\alpha\mathbb{Z}_5 \times 13\mathbb{Z}_{13} \times \cdots \times p_h\mathbb{Z}_{p_h} \times \mathbb{Z}_{p_{h+1}} \times \cdots \times \mathbb{Z}_{p_k}.\]
We prove the following variation of Proposition 3.5 in \cite{kravitz2025}. 
\begin{prop}\label{prop:f_h}
    Let $n, h, k, \alpha$ be positive integers with $h \le k$. Then $F(n, h, k, \alpha)$ is a fundamental domain of $X/(\theta_5^{-n}A)$. 
\end{prop}
\begin{proof}
    Since $X = (\mathbb{C} \times \mathbb{Q}_5 \times \cdots \times \mathbb{Q}_{p_k})/A$, we deduce that
    \[X/(\theta_5^{-n}A) = (\mathbb{C} \times \mathbb{Q}_5 \times \cdots \times \mathbb{Q}_{p_k})/(\theta_5^{-n}A) = (\mathbb{C} \times \mathbb{Q}_5 \times \cdots \times \mathbb{Q}_{p_k})/(P_5^{-n}A).\]
    By Lemma \ref{lmm:intertwine}, we have that $A$ has index $5^n$ as a subgroup of $\theta_5^{-n}A$, which acts freely on $\mathbb{C} \times \mathbb{Q}_5 \times \cdots \times \mathbb{Q}_{p_k}$. Thus,
    \[\text{vol}(X/(\theta_5^{-n}A)) = 5^{-n} \text{vol}(X) = 5^{-n} = \text{vol}(F(n,h,k,\alpha)).\]

    It suffices to show that every element of $\mathbb{C} \times \mathbb{Q}_5 \times \cdots \times \mathbb{Q}_{p_k}$ has a representative modulo $\theta_5^{-n}A$ within $F$. Indeed, by the Chinese remainder theorem (or the Strong Approximation Theorem, see Chapter 10 of Cassels \cite{cassels1986}), we can subtract an element of $A$ to restrict our attention to
    \[(z, a_5, \dots, a_{p_k}) \in \mathbb{C} \times \mathbb{Q}_5 \times \cdots \times \mathbb{Q}_{p_k}, \quad a_5 \in 5^\alpha \mathbb{Z}_5, \quad \forall 2 \le j \le h, a_{p_j} \in p_j\mathbb{Z}_{p_j}, \quad \forall h < j \le k, a_{p_j} \in \mathbb{Z}_{p_j}.\]
    Now since
    \[P_5^{-n}A \cap 5^\alpha \mathbb{Z}_5 \cap \left(\bigcap_{j = 2}^h p_j\mathbb{Z}_{p_j}\right) \cap \left(\bigcap_{j = h+1}^k \mathbb{Z}_{p_j}\right) = \left(P_5^{-n}\overline{P_5}^\alpha \prod_{j = 2}^h \overline{P_{p_j}}\right)\mathbb{Z}[i],\]
    we can further subtract an element of the right-hand lattice to preserve all the $p$-adic valuations, yet shift $z$ to any fundamental domain, say,
    \[\left(P_5^{-n}\overline{P_5}^\alpha \prod_{j = 2}^h \overline{P_{p_j}}\right) ([-1/2,1/2) \times [-1/2,1/2)i)).\qedhere\]
\end{proof}

Let $\mathcal{P}(n, h, k, \alpha)$ denote the partition of $X$ given by the $(\theta_5^{-n}A)$-orbits of $F(n, h, k, \alpha)$. Additionally, setting $n = 0$ in Proposition \ref{prop:f_h} tells us that
\[F_{h, k, \alpha} := \left(\overline{P_5}^\alpha \prod_{j = 2}^h \overline{P_{p_j}}\right) ([-1/2,1/2) \times [-1/2,1/2)i)) \times 5^\alpha\mathbb{Z}_5 \times 13\mathbb{Z}_{13} \times \cdots \times p_h\mathbb{Z}_{p_h} \times \mathbb{Z}_{p_{h+1}} \times \cdots \times \mathbb{Z}_{p_k}\]
is a fundamental domain of $X$ itself. 

Since parts of $\mathcal{P}(n, h, k, \alpha)$ look like, for $w \in \mathbb{Z}[i]$, 
\[\left(P_5^{-n}\overline{P_5}^\alpha \prod_{j = 2}^h \overline{P_{p_j}}\right) (w + [-1/2,1/2) \times [-1/2,1/2)i)) \times 5^\alpha\mathbb{Z}_5 \times 13\mathbb{Z}_{13} \times \cdots \times p_h\mathbb{Z}_{p_h} \times \mathbb{Z}_{p_{h+1}} \times \cdots \times \mathbb{Z}_{p_k}.\]
We note that $\mathcal{P}(n, h, k, \alpha)$ gives a direct partition of $F_{h, k, \alpha}$. Furthermore, the difference of any two parts, indexed by $w, w'$, is the union of the four parts indexed by
\[w - w', \quad w - w' - 1, \quad w - w' - i, \quad w - w' -1 - i.\]

\subsubsection{The group structure of $(\theta_5^{-n}A)/A$.} \label{subsec:theta5} In Section \ref{subsec:entropy}, the group $(\theta_5^{-n}A)/A$ plays a crucial role as a discrete approximation to our partition. We develop some background on the group that allows us to adapt the work of Bourgain, Lindenstrauss, Michel and Venkatesh \cite{bourgain2009}. 

We consider the following lemmas, analogous to Section 3.1 of \cite{kravitz2025}. The first allows us to intertwine $(\theta_5^{-n}A)/A$ with $\mathbb{Z}/5^n\mathbb{Z}$, with statement identical to Lemma 3.1 of \cite{kravitz2025}. 

\begin{lmm}\label{lmm:intertwine}
    Let $n$ be a non-negative integer. Then
    \[(\theta_5^{-n}A)/A = (P_5^{-n}\mathbb{Z}[i])/A \cong \mathbb{Z}[i]/P_5^n\mathbb{Z}[i] \cong \mathbb{Z}/5^n\mathbb{Z},\]
    with the final isomorphism chosen to intertwine the $\times P_5^\ell$ and $\times 5^\ell$ maps for all $0 \le \ell \le n$. 
\end{lmm}
\begin{proof}
    For the first equality, it suffices to show that $P_5^{-n}\mathbb{Z}[i] + A = \theta_5^{-n}A + A$. We have
    \[P_5^{-n}\mathbb{Z}[i] = \theta_5^{-n}\overline{P_5}^n\mathbb{Z}[i] \subseteq \theta_5^{-n}A.\]
    This gives one inclusion. For the other, write $a \in A$ as
    \[a = \prod_{j = 1}^k\overline{P_{p_j}}^{\alpha_j} b\]
    for $b \in \mathbb{Z}[i]$ coprime to each $\overline{P_{p_j}}$, $1 \le j \le k$, and $\alpha_j \in \mathbb{Z}$. Now consider adding
    \[c = \overline{P_5}^n \prod_{j = 1}^k\overline{P_{p_j}}^{\alpha_j}d\]
    for some $d \in \mathbb{Z}[i]$. Since $P_5$ is coprime to each of $\overline{P_{p_j}}$, we can choose $d$ such that
    \[\overline{P_5}^n \prod_{j = 1}^k\overline{P_{p_j}}^{\alpha_j} (b + P_5^nd) \in \mathbb{Z}[i].\]
    Dividing through by $P_5^n$, we obtain
    \[\theta_5^na + c \in P_5^{-n}\mathbb{Z}[i],\]
    which gives the reverse inclusion $\theta_5^{-n}A + A \subseteq P_5^{-n}\mathbb{Z}[i] + A$.

    The first isomorphism is given by, for $u \in \mathbb{Z}[i]$,
    \[P_5^{-n}u + A \leftrightarrow u + P_5^n\mathbb{Z}[i].\]
    Indeed, $P_5^{-n}(u-u') \in A$ if and only if $u - u' \in P_5^nA$. Since $u, u' \in \mathbb{Z}[i]$, this is equivalent to $u - u' \in P_5^n\mathbb{Z}[i]$. For the second isomorphism, consider
    \[m + \mathbb{Z}[i] \leftrightarrow m + 5^n\mathbb{Z}\]
    for $0 \le m \le 5^n-1$. Since
    \[\{m + \mathbb{Z}[i] : 0 \le m \le 5^n-1\}\]
    form a reduced residue system modulo $P_5^n$ over $\mathbb{Z}[i]$, this yields an isomorphism that intertwines $\times P_5^\ell$ and $\times 5^\ell$ for $0 \le \ell \le n$. 
\end{proof}

We also wish to consider Fourier analysis on $(\theta_5^{-n}A)/A$. Thus, in the next lemma, we describe the characters on this group. This generalises Lemma 3.3 of \cite{kravitz2025}. 

\begin{lmm}\label{lmm:dual}
    The dual of $(\theta_5^{-n}A)/A$ is isomorphic to $A/(\theta_5^nA)$ via the map $\Psi$ defined as follows. For $\xi \in A/(\theta_5^nA)$ and $x \in (\theta_5^{-n}A)/A$, we have
    \[\Psi(\xi) : x \mapsto -\Re(\xi x) + \sum_{j = 1}^k \{\iota_{p_j}(\xi \cdot x)\}_{p_j} \pmod{1}.\]
\end{lmm}
\begin{proof}
    Note that we have, for $z \in \mathbb{Z}[i]$,
    \[-\Re(z) + \sum_{j = 1}^k \{\iota_{p_j}(z)\}_{p_j} \in \mathbb{Z}\]
    if and only if $z \in A$. Thus, our map is well-defined since we can shift $x$ by arbitrary elements of $A$, and it is injective since evaluation at $x = \theta_5^{-n}$ would give $\xi \in \theta_5^nA$. $\Psi$ is also clearly additive.

    Note that both groups have cardinality equal to $5^n$ (as the dual group of a finite abelian group has the same cardinality). Since $\Psi$ is an injective homomorphism on finite groups of the same cardinality, it is an isomorphism. 
\end{proof}

\subsubsection{Hecke characters.} \label{subsec:hecke} We briefly introduce the relevant theory of Hecke characters for $\mathbb{Q}(i)$. For further exposition into Größencharaktere over general number fields, we refer to the works of Iwaniec, Kowalski \cite{iwaniec2004} and Green, Sawhney \cite{green2026}. We apply this theory in Section \ref{subsec:quasirh}. 

To distinguish between types of Gaussian primes, let $\pi$ always be such that $N(\pi) = p$ prime, and let $\mathfrak{p}$ denote a general Gaussian prime. We only consider Dirichlet characters, which we refer to as finite Hecke characters to avoid confusion, as there are also angular (or infinite) characters. 

For $d \in \mathbb{Z}[i]$, we consider characters modulo $d$,
\[\chi : (\mathbb{Z}[i]/d\mathbb{Z}[i])^\times \to \mathbb{C}^\times,\]
with $\chi(\alpha) = 0$ if $\gcd(\alpha, d) \neq 1$. The conductor of $\chi$ is the minimal $d$ (in norm) such that this holds. 

We also define the Hecke $L$-function
\[L(s, \chi) := \sum_{0 \neq \alpha \in \mathbb{Z}[i]} \frac{\chi(\alpha)}{N(\alpha)^s},\]
as well as the analogue of the Chebyshev functions
\[\psi(x, \chi) := \sum_{N(\alpha) \le x} \Lambda(\alpha) \chi(\alpha), \quad \vartheta(x, \chi) := \sum_{N(\mathfrak{p}) \le x} \Lambda(\mathfrak{p}) \chi(\mathfrak{p}),\]
where:
\[\Lambda(\alpha) := \begin{cases}
    \log N(\mathfrak{p}) & \alpha = \mathfrak{p}^k \text{ for prime $\mathfrak{p}$, integer $k \ge 0$} \\
    0 & \text{otherwise.}
\end{cases}\]

Theorem 2.9.2 of Schlackow \cite{schlackow2010} gives us the following explicit formula.

\begin{lmm}\label{lmm:explicit}
    Let $\chi$ be a non-principal character with conductor $d$ and $2 \le T \le x$. Then
    \[\psi(x, \chi) = -\frac{x^{\beta_1}}{\beta_1} - \sum_{\substack{|\Im(\rho)| \le T\\\rho \neq \beta_1, 1-\beta_1}} \frac{x^\rho}{\rho} + O\left(\frac{x}{T} \log(|d|x)^2 + x^{1/4} \log x\right),\]
    where:
    \begin{outline}[enumerate]
        \1 The term $-\frac{x^{\beta_1}}{\beta_1}$ is excluded unless $\chi$ is a real character for which $L(s, \chi)$ has a real Siegel zero $\beta_1$ such that $\beta_1 > 1 - \frac{F}{\log |d|}$ for some positive constant $F$.
        \1 The sum is over all non-trivial zeros $\rho$ of $L(s,\chi)$ which are not Siegel with $|\Im(\rho)| \le T$. 
    \end{outline}
\end{lmm}

Very recently, there have been developments for the zero-free strips of Hecke $L$-functions. The initial works generated by OpenAI \cite{openai2026primitive, openai2026twelve, openai2026quasirh} demonstrated there are no zeros of Dirichlet $L$-functions and finite Hecke $L$-functions over $\mathbb{Q}(\sqrt{-3})$ in the half-plane $\Re(s) > 7/8$. Subsequently Cureton and Thorner used AxiomProver to prove and Lean-formalise a zero-free half-plane for finite Hecke $L$-functions over arbitrary number fields in \cite{thirtytwo}. 

\begin{thm}\label{thm:thirtytwo}
    Let $\chi$ be a finite Hecke character defined over a number field $K$. Then its Hecke $L$-function $L(s, \chi)$ does not vanish when $\Re(s) > 31/32$, wherever it is holomorphic. 
\end{thm}

This has been improved to $\Re(s) > 7/8$ by Hu using AI assistance \cite{seveneigths}. We remark that any zero-free strip transfers to quasipolynomial pyjama bounds, with exponent depending on the threshold. Hence, Lemma \ref{thm:thirtytwo} suffices for our purposes. In particular, there are no Siegel zeros, and Lemma \ref{lmm:explicit} converts into the following unconditional statement.

\begin{cor}\label{cor:explicit2}
    Let $\chi$ be a non-principal character with conductor $d$ and $2 \le T \le x$. Then
    \[\psi(x, \chi) = -\sum_{|\Im(\rho)| \le T} \frac{x^\rho}{\rho} + O\left(\frac{x}{T} \log(|d|x)^2 + x^{1/4} \log x\right),\]
    where the sum is over all non-trivial zeros $\rho$ of $L(s, \chi)$ with $|\Im(\rho)| \le T$. These zeros satisfy $\Re(\rho) \le 31/32$. 
\end{cor}

\subsection{Initial pigeonhole} \label{sec:pigeon}

We wish to apply an initial pigeonhole step which is more height-friendly than \cite{kravitz2025}. In particular, Kravitz and Leng consider a set of $N_0^2$ rotations which yields a difference set having height at most $\exp(100N_0)$. An improvement can be achieved if we work with continuous characters, where the group action can be expanded to encompass all of $\mathbb{Q}(i)$.

\begin{prop}\label{prop:pigeon}
    Let $N_0$ be a positive integer. Then there exists $\Theta_0 \subseteq \mathbb{Q}(i) \cap S^1$ of cardinality $N_0$ such that the following holds. Let $p := \chi_z$ be a continuous character. Then there exists $r \in \Theta_0 - \Theta_0$ of height at most $N_0^2+1$ such that
    \[\dist(rp, 0) \le (k+2)N_0^{-1/(k+2)}.\]
\end{prop}
\begin{proof}
    Note that our initial group action, 
    \[(\theta \cdot \chi)(z) := \chi(\theta z)\]
    only makes sense if $\theta \in A$. Nonetheless, if we restrict to the subset of continuous characters, we can in turn expand our set of rotations. Consider the action on $j_\mathbb{C}(\mathbb{C})$ given by, for $\phi \in \mathbb{Q}(i)$,
    \[\phi \circ \chi_z := \chi_{\phi z}.\]
    This coincides with the usual one on $A$. This is well-defined since each continuous character has a unique representation as $\chi_z$. Indeed, if $\chi_z = \chi_{z'}$, then $\Re((z - z')w) = 0$ for all $w \in \mathbb{C}$. Setting $w = \overline{z - z'}$ gives $|z - z'|^2 = \Re(|z - z'|^2) = 0$, hence $z = z'$. 

    We can then define
    \[\Theta_0 := \{1\} \cup \left\{\frac{m+i}{m-i} : 1 \le m \le N_0-1\right\}.\]
    Note that any two elements in $\Theta_0$ are distinct. Indeed, 1 is the only real element, and otherwise
    \begin{align*}
        \frac{m+i}{m-i} - \frac{m'+i}{m'-i} &= \frac{(m+i)(m'-i) - (m-i)(m'+i)}{(m-i)(m'-i)} \\
        &= \frac{2i(m' - m)}{(m-i)(m'-i)},
    \end{align*}
    hence if the left-hand side vanishes, then $m = m'$. Additionally, we find that $\Theta_0 - \Theta_0$ has maximum height
    \[\max\{2|m' - m|, |(m-i)(m'-i)|\} \le N_0^2 + 1.\]

    By the pigeonhole principle within $F_0$, we can find distinct $\phi_1, \phi_2 \in \Theta_0$ such that
    \[\dist((\phi_1 - \phi_2)p, 0) = \dist(\phi_1p, \phi_2p) \le (k+2)N_0^{-1/(k+2)}.\]
    Let $r := \phi_1 - \phi_2$. Thus, there exists a representative $rp = j(z, a_5, \dots, a_{p_k})$ where
    \[|z|_\mathbb{C} + \sum_{j = 1}^k |a_{p_j}|_{p_j} \le (k+2)N_0^{-1/(k+2)}.\qedhere\]
\end{proof}

\subsection{Major arcs} \label{sec:major}

Recall that the major-arc case represents when $rp$ (in the language of Proposition \ref{prop:pigeon}) lies extremely close to $\overline{B}_1^\mathbb{C}(0)$. Intuitively, we wish to clear denominators to conclude that $p$ lies close to a torsion point of small order. However, one must also expand the complex ball to allow this. We continue similarly to \cite{kravitz2025} by considering casework on the size of the complex ball, proving that a further orbit of $p$ intersects $E^*(\eps)$ if the radius is too large.

\subsubsection{Baker's Theorem.} \label{subsec:baker} We use Baker's Theorem \cite{bakerI1967, bakerII1967, bakerIII1967} to provide quantitative approximation results. For $\beta$ an algebraic number, let $H(\beta)$ denote its height. 

\begin{thm}\label{thm:baker}
    Suppose $\lambda_1, \dots, \lambda_n \in \mathbb{C}$ are linearly independent over the rationals, with $e^{\lambda_i} \in \overline{\mathbb{Q}}$ for $1 \le i \le n$. Then there exists a constant $c > 0$ depending on $\lambda_1, \dots, \lambda_n$, such that for any algebraic numbers $\beta_1, \dots, \beta_n$ not all zero, we have
    \[\left|\sum_{i = 1}^n \lambda_i \beta_i\right| \ge c\max\{H(\beta_1), \dots, H(\beta_n)\}^{-c}.\]
\end{thm}

We then extract the following corollary, identical to Corollary 4.2 of \cite{kravitz2025}. Let $\alpha_5 \in [-1/2, 1/2)$ such that $e^{2\pi i\alpha_5} = \theta_5$. 

\begin{cor}\label{cor:baker}
    There exists a constant $c' > 1$ such that for every non-zero rational $q$ of height at most $H$, we have
    \[||\alpha_5 - q||_{\mathbb{R}/\mathbb{Z}} \ge c'H^{-c'}.\]
\end{cor}

We provide an alternative direct proof of the following ``quantitative approximation lemma", Lemma 4.3 of \cite{kravitz2025}. This tells us that we can obtain $\eps$-density within $S^1$ using only $\text{poly}(\eps^{-1})$ powers of $\theta_5$. We use this both in the major arcs argument, as well as in the final conversion to the plane in Section \ref{sec:conversion}. 

\begin{lmm}\label{lmm:quantapprox}
    There exists a constant $c > 0$ such that the following holds. If $K, \ell \in \mathbb{N}$ with $1 \le \ell \le K^c$, then
    \[\{m\ell \alpha_5 : 0 \le m \le K\}\]
    is $K^{-c}$-dense in $\mathbb{R}/\mathbb{Z}$. 
\end{lmm}
\begin{proof}
    This immediately holds for $K = 1$. Henceforth, assume $K \ge 2$. 

    Let $c_1, c_2 > 0$ to be determined later, with $1 \le \ell \le K^{c_2}$. Consider $\{m\ell\alpha_5 : 0 \le m \le K^{c_1}\}$. By the pigeonhole principle, we can find distinct $0 \le m, m' \le K^{c_1}$ with:
    \[||(m - m')\ell\alpha_5||_{\mathbb{R}/\mathbb{Z}} \le K^{-c_1}.\]
    Let $d = m - m'$, so $|d| \le 2K^{c_1}$. For some integer $a$, we have
    \[|d\ell \alpha_5 - a| \le K^{-c_1}.\]

    Note that by the triangle inequality, $|a| \le |d\ell\alpha_5| + K^{-c_1} \le |d\ell|$. Thus, by Corollary \ref{cor:baker}, there exists a constant $c' > 1$ such that
    \begin{align*}
        |d\ell\alpha_5 - a| = \left|\alpha_5 - \frac{a}{d\ell}\right| \cdot |d\ell| &\ge |d\ell| c'\max\{|a|, |d\ell|\}^{-c'} \\
        &\ge c'|d\ell|^{1-c'} \\
        &\ge c'2^{1-c'}K^{(1-c')(c_1+c_2)}.
    \end{align*}
    In particular, a further orbit of $\{\tilde{m} \cdot |d|\ell\alpha_5 : 0 \le \tilde{m} \le (c')^{-1} 2^{c'-1}K^{(c'-1)(c_1+c_2)}\}$ is $K^{-c_1}$-dense within $\mathbb{R}/\mathbb{Z}$. In total,
    \[\{m\ell\alpha_5 : 0 \le m \le (c')^{-1}2^{c'-1} K^{(c'-1)(c_1+c_2) + c_1}\}\]
    is $K^{-c_1}$-dense. Picking $c_1, c_2$ sufficiently small and setting $c := \min\{c_1, c_2\}$ proves the result. 
\end{proof}

\subsubsection{Geometric lemmas.} \label{subsec:geometric} As in Section 4.2 of \cite{kravitz2025}, we consider three geometric lemmas. Only the first lemma on clearing denominators, analogous to Lemma 4.4 of \cite{kravitz2025}, requires us to spell out the proof again; the others follow almost identically. 

\begin{lmm}\label{lmm:geo1}
    Let $\varpi > 0$ and $H > 1$. Suppose $p \in X$ satisfies $\dist(rp, \overline{B}_1^\mathbb{C}(0)) \le \varpi$ for some non-zero Gaussian rational $r \in \mathbb{Q}(i)$ of height at most $H$. Then there is a torsion point $q \in X$ of order at most $H^2$ such that
    \[\dist(p, \overline{B}_H^\mathbb{C}(q)) \le (k+1)H^2\varpi.\]
    \begin{proof}
        We can write $rp = j(z, a_5, \dots, a_{p_k})$ with $|a_{p_j}|_{p_j} \le \varpi$ for each $1 \le j \le k$ and $|z|_{\mathbb{C}} \le 1 + \varpi$. For all $1 \le j \le k$, we have that $|\iota_{p_j}(r)|_{p_j} \ge H^{-2}$. Thus, we can set $a_{p_j}' := a/\iota_{p_j}(r)$ with $|a_{p_j}'|_{p_j} = |a_{p_j}|_{p_j} \cdot |\iota_{p_j}(r)|_{p_j}^{-1} \le H^2\varpi$. Similarly, let $z' := z/r$ with $|z'|_\mathbb{C} \le H|z|_\mathbb{C} \le H(1 + \varpi)$. 

        Now, we have that
        \[r(p - j(z', a_5', \dots, a_{p_k}')) = j(z - rz', a_5 - ra_5', \dots, a_{p_k} - ra_{p_k}) = 0,\]
        hence $q := p - j(z', a_5', \dots, a_{p_k}')$ is a torsion point of order $|r|^2 \le H^2$. Our bounds on $|a_{p_j}'|_{p_j}$ and $|z'|_\mathbb{C}$ above imply that $p$ lies within distance
        \[H\varpi + \sum_{j = 1}^k H^2\varpi \le (k+1)H^2\varpi\]
        of $\overline{B}_H^\mathbb{C}(q)$, as desired.
    \end{proof}
\end{lmm}

The next lemma tells us we can multiply up torsion points until they are invariant under some power of $\theta_5$. It is proven identically to Lemma 4.5 of \cite{kravitz2025} and is reminiscent of Corollary 6.5 of \cite{manners2013}. 

\begin{lmm}\label{lmm:geo2}
    Let $m \in \mathbb{N}$. If $q \in X$ is a torsion point of order $m$, then there are integers $0 \le s_0 < m^2$ and $1 \le s_1 \le m^2$ such that $q' := \theta_5^{s_0}q$ satisfies $\theta_5^{s_1}q' = q'$. 
\end{lmm}

Finally, we note that any sufficiently large complex circle intersects the diagonal set $E^*(0) := \{\chi \in X : \chi(1) = 0\}$. This is Lemma 4.6 of \cite{kravitz2025}. 

\begin{lmm}\label{lmm:geo3}
    For every radius $r \ge 1/2$ and $p \in X$, the complex circle $S_r^\mathbb{C}(p)$ intersects $E^*(0)$. 
\end{lmm}

\subsubsection{Major arc dichotomy.} \label{subsec:major} Now we combine all the ingredients together. Given a point $rp$ extremely close to $\overline{B}_1^\mathbb{C}(0)$, we apply Lemma \ref{lmm:geo1} to deduce that $p$ is very close to a complex ball, centred at a torsion point of low order. We are done if the radius is at most 1/2. Otherwise, we use Lemmas \ref{lmm:geo2} and \ref{lmm:geo3} to show that a $\theta_5$-orbit of $p$ intersects $E^*(\eps)$. 

\begin{prop}\label{prop:major}
    There exists a constant $C > 0$ such that the following holds for all choices of $0 < \eps_1, \delta < 1/10$ and $H \in \mathbb{N}$. Set $K = (H\eps_1^{-1})^C$ and suppose that
    \[\dist(rp, \overline{B}_1^\mathbb{C}(0)) \le \varpi := \min\left\{\frac{\delta}{(k+1)H^2}, \exp(-10K)\right\}\]
    for some non-zero $r \in \mathbb{Q}(i)$ of height at most $H$. Then either: we have $\theta_5^sp \in E^*(\eps_1)$ for some $0 \le s \le K$; or there is a torsion point $q \in X$ of order at most $H^2$ such that $\dist(p, \overline{B}_{1/2}^\mathbb{C}(q)) \le \delta$.
\end{prop}
\begin{proof}
    Given Lemmas \ref{lmm:quantapprox}, \ref{lmm:geo1}, \ref{lmm:geo2}, \ref{lmm:geo3}, we can proceed almost identically to the proof of Proposition 4.7 in \cite{kravitz2025}, with the only difference being the factor of $(k+1)H^2$ instead of (re-phrased in our language) $2H^2 + H$ to reflect our use of $k$ primes instead of two. 
\end{proof}

\subsection{Minor arcs} \label{sec:minor}

Our minor arc case is split into two mixing stages, as in Section 5 of the Kravitz-Leng framework \cite{kravitz2025}. Since $rp$ is not extremely close to $\overline{B}_1^\mathbb{C}(0)$, by construction of our metric, some $p_j$-adic component must be slightly large. For convenience, we select the prime 5, but we will check in Section \ref{subsec:minor} that the general case follows identically and without affecting the shape of the results. 

In Section \ref{subsec:quasirh}, we prove our desired number theoretic input about generating multiplicative groups modulo powers of $\overline{P_5}$ using subset products of the other $\theta_{p_j}$. This serves to replace the ``Lifting the Exponent" lemma used in \cite{kravitz2025}, \cite{manners2013}. Using subset products allows us to work with noticeably smaller distances. 

In Section \ref{subsec:manners}, the first mixing stage, we apply the above input as a density statement on discrete scales. Indeed, we ``scale down" using powers of 5, then mix with the other primes to obtain a notion of density which is much weaker than intersecting every $\eta$-ball. Yet this notion both observes the pyjama stripe and gives leeway for a second subset of mixing primes. 

In Section \ref{subsec:entropy}, the second mixing stage, we adapt the machinery of Bourgain, Lindenstrauss, Michel, Venkatesh \cite{bourgain2009} using our number theoretic input and a bump function catered to $E^*(\eps)$. Here we take a smaller power of 5 as a scale, which requires a smaller set of mixing primes. 

It would be possible to run our arguments with two disjoint sets of primes instead, one for each mixing stage. Quantitatively, our primes will be roughly the same size, and so a similar $p$-adic analysis follows through. Thus, the primary statement to adapt would be generation using subset products of a non-initial segment of small primes. Yet by running Lemma \ref{lmm:character-sum} twice with different parameters and then taking the difference, we obtain the power-saving result we require. 

\subsubsection{Subset products of $\theta_p$.} \label{subsec:quasirh} Our proof roughly follows the structure of Klurman, Shparlinski, Teräväinen \cite{klurman2026}, which notably works within the usual integers $\mathbb{Z}$ and does not assume any zero-free strip. Instead, they appeal to the grand density theorem (Theorem 10.4 of \cite{iwaniec2004}), which yields an exceptional set of excluded moduli. 

To demonstrate the generation of all coprime residues via subset products, one expands out the relevant indicator functions using character orthogonality, then it suffices to bound a character sum ranging over primes. Here, Theorem \ref{thm:thirtytwo} bypasses many steps of the proof, as one would previously expect when assuming the Generalised Riemann Hypothesis. 

\begin{lmm}\label{lmm:character-sum}
    Let $Q$ be a non-zero Gaussian integer coprime to $\overline{Q}$. Let $\chi$ be a non-principal finite Hecke character modulo $Q$. Then
    \[\left|\sum_{N(\pi) \le x} \chi(\pi/\overline{\pi})\right| \ll x^{31/32} (\log (N(Q)x))^2.\]
\end{lmm}
\begin{proof}
    Note that
    \[\chi(\pi/\overline{\pi}) = \chi(\pi) \overline{\chi}(\overline{\pi}) = \chi(\pi)\tilde{\chi}(\pi),\]
    where $\tilde{\chi}(\pi) := \overline{\chi}(\overline{\pi})$. This is a bona fide finite Hecke character modulo $\overline{Q}$, indeed we have
    \[\tilde{\chi}(\alpha + m\overline{Q}) = \overline{\chi}(\overline{\alpha} + \overline{m}Q) = \overline{\chi}(\overline{\alpha}), \quad \tilde{\chi}(\alpha_1)\tilde{\chi}(\alpha_2) = \overline{\chi}(\overline{\alpha_1})\overline{\chi}(\overline{\alpha_2}) = \tilde{\chi}(\alpha_1\alpha_2).\]
    Thus, $\chi(\pi/\overline{\pi}) = \gamma(\pi)$ for $\gamma$ a finite Hecke character modulo $Q\overline{Q} = N(Q)$. Furthermore, we claim $\gamma$ is non-principal. Indeed, suppose for the sake of contradiction that $\gamma$ is principal. Then modulo $N(Q)$, we have that $\overline{\chi}$ and $\tilde{\chi}$ are the same character. Thus, they are induced by the same character modulo $\gcd(Q, \overline{Q}) = 1$, i.e. principal themselves. Yet we assumed $\chi$ was not principal. 

    Now we consider $\psi(x, \gamma)$ and $\vartheta(x, \gamma)$. Picking $T = x^{1/32}$ in Corollary \ref{cor:explicit2}, we obtain
    \[\psi(x,\gamma) = -\sum_{|\Im(\rho)| \le x^{1/32}} \frac{x^\rho}{\rho} + O(x^{31/32} \log(N(Q)x)^2 + x^{1/4}\log x).\]
    Equations (2.32) and (2.33) in Schlackow \cite{schlackow2010} give us the bound of
    \[\sum_{|\Im(\rho)| \le x^{1/32}} \frac{1}{|\rho|} \ll \log(N(Q)x)^2.\]
    Since $|x^\rho| \le x^{31/32}$, we combine our estimates to get
    \[|\psi(x, \gamma)| \ll x^{31/32} \log(N(Q)x)^2.\]
    
    Applying partial summation, we find that
    \[\sum_{N(\mathfrak{p}) \le x} \gamma(\mathfrak{p}) = \frac{\vartheta(x, \gamma)}{\log x} + \int_2^x \frac{\vartheta(t, \gamma)}{t(\log t)^2} dt \ll x^{31/32} \log(N(Q)x)^2.\]
    Now note that if $\mathfrak{p} = p \equiv3 \pmod{4}$ is a rational prime, then $N(p) = p^2$. Thus,
    \[\left|\sum_{N(p) = p^2 \le x} \gamma(p)\right| \le \sum_{p \le x^{1/2}} 1 = O(x^{1/2}).\]
    Altogether, we obtain
    \[\left|\sum_{N(\pi) \le x} \chi(\pi/\overline{\pi})\right| = \left|\sum_{N(\pi) \le x} \gamma(\pi)\right| \ll x^{31/32} \log(N(Q)x)^2.\qedhere\]
\end{proof}

This leads to a discrepancy inequality, which we use to bound the non-principal parts of the character orthogonality decomposition. We construct our set of subset products
\[M_k := \left\{\prod_{j = 2}^k \theta_{p_j}^{e_j} : e_j \in \{0,1\}\right\}.\]

\begin{thm}\label{thm:subset-products}
    Let $Q$ be a non-zero Gaussian integer coprime to $\overline{Q}$. Suppose it is also coprime to $p_2, \dots, p_k$, where $p_j$ is the $j$th rational prime satisfying $p_j \equiv 1 \pmod{4}$. Furthermore, suppose that $N(Q) \le \exp(k^{1/65})$. Then we have
    \[\left\{\prod_{j = 2}^k \theta_{p_j}^{e_j} \pmod{Q} : e_j \in \{0,1\}\right\} = (\mathbb{Z}[i]/Q\mathbb{Z}[i])^\times.\]
\end{thm}
\begin{proof}
    For a non-principal character $\chi$ mod $Q$, consider the quantity
    \[D(\chi) = \sum_{j = 2}^k (1 - \Re \chi(\theta_{p_j})).\]
    Note that we have
    \[\Re\chi(\theta_p) = \frac12 (\chi(\theta_p) + \overline{\chi(\theta_p)}) = \frac12 (\chi(\theta_p) + \chi(\overline{\theta_p})).\]
    Hence, we deduce from Lemma \ref{lmm:character-sum},
    \begin{align*}
        D(\chi) &= (k-1) - \frac12 \sum_{j = 2}^k (\chi(\theta_{p_j}) + \chi(\overline{\theta_{p_j}})) \\
        &= k-1 - \frac12 \sum_{5 < N(\pi) \le p_k} \chi(\pi/\overline{\pi}) \\
        &\ge k-1 + O(k^{31/32} \log(N(Q)k)^2).
    \end{align*}
    Thus, for $N(Q) \le \exp(k^{1/65})$, we have $D(\chi) \ge k/2$, say. 

    For $a \in (\mathbb{Z}[i]/Q\mathbb{Z}[i])^\times$, define the indicator
    \[T_Q(a) := \sum_{e_j \in \{0,1\}} 1_{\prod_{j = 2}^k \theta_{p_j}^{e_j} \equiv a \pmod{Q}}.\]
    Let $\varphi(Q) = |(\mathbb{Z}[i]/Q\mathbb{Z}[i])^\times|$. By orthogonality of characters, we find that
    \begin{align*}
        T_Q(a) &= \frac{1}{\varphi(Q)} \sum_{\text{$\chi$ mod $Q$}} \overline{\chi}(a) \prod_{j = 2}^k (1 + \chi(\theta_{p_j})) \\
        &= \frac{2^k}{\varphi(Q)} + \frac{1}{\varphi(Q)} \sum_{\substack{\text{$\chi$ mod $Q$}\\\chi \neq \chi_0}} \overline{\chi}(a) \prod_{j = 2}^k (1 + \chi(\theta_{p_j})).
    \end{align*}
    
    To bound the non-principal contribution, we obtain from the arithmetic-geometric mean inequality that
    \begin{align*}
        \left|\prod_{j = 2}^k (1 + \chi(\theta_{p_j}))\right| &\le \left(\frac{1}{k-1} \sum_{j = 2}^k |1 + \chi(\theta_{p_j})|^2\right)^{(k-1)/2} \\
        &= \left(\frac{1}{k-1} \sum_{j =2}^k (2 + 2\Re \chi(\theta_{p_j}))\right)^{(k-1)/2} \\
        &= \left(4 - \frac{2D(\chi)}{k-1}\right)^{(k-1)/2} \\
        &= 2^{k-1} \left(1 - \frac{D(\chi)}{2(k-1)}\right)^{(k-1)/2} \\
        &\le 2^{k-1} \exp\left(-\frac{D(\chi)}{4}\right) \\
        &\le 2^{k-1} \exp(-k/8).
    \end{align*}
    Since $\exp(k/8) \ge N(Q) \ge \varphi(Q)$, the principal contribution dominates and $T_Q(a) >0$ for every $a \in (\mathbb{Z}[i]/Q\mathbb{Z}[i])^\times$. 
\end{proof}

Note that an identical statement holds if we single out any other prime $p_j$. 

\subsubsection{Manners's denseness criterion.} \label{subsec:manners} This is the first mixing stage. As described in Section \ref{subsec:improvements}, we weaken the notion of density within our growing space $X$ in light of the second mixing stage. Instead, we focus on intersecting parts of our partition $\mathcal{P}(n, h, k, \alpha)$, which are given by translates of $F(n, h, k, \alpha)$. Recall their definitions from Section \ref{subsec:fundamental}.  

We follow the general proof outline of Theorem 5.1 of \cite{kravitz2025}, obtained by analysing Manners's treatment of ``non-Archimedean limit points" in Lemma 7.1 of \cite{manners2013}. Let $\Xi_h := \prod_{j = 2}^h p_j$. 

\begin{thm}\label{thm:manners-denseness}
    Let $0 \le \eta, \varpi \le 1/10$ and $2 \le h \le k$ be integers. Suppose $x \in X$ can be represented as $x = j(z, a_5, a_{13}, \dots, a_{p_k})$ where
    \[|z|_\mathbb{C} \le \frac{\eta}{2}, \quad \frac{\varpi}{k} \le |a_5|_5 \le \frac{\eta}{2}, \quad \forall 2 \le j \le h, |a_{p_j}|_{p_j} \le \frac{1}{p_j^2}, \quad \forall h < j \le k, |a_{p_j}|_{p_j} \le \frac{1}{p_j}.\]
    Then for integers $\alpha, n$ satisfying
    \[\alpha \le k^{1/80} + \log_5(\eta^2\Xi_h^{-1}/2^8),\quad \eta \le (5^{\alpha - n}\Xi_h)^{1/2}/2,\]
    the orbit given by
    \[\{\theta_5^smx : m \in M_k, 0 \le s \le S_1\}, \quad S_1 = \log_5(k\varpi^{-1}\eta^{-2}\Xi_h/2^8)\]
    intersects every set of the form $x' + F(n, h, k, \alpha)$ for $x' \in X$. 
\end{thm}
\begin{proof}
    For $\alpha' := \ceil{\log_5(\eta^2\Xi_h^{-1}/2^8)}$, we can observe $x'$ in the fundamental domain $F_{h,k,\alpha'}$ given in Section \ref{subsec:fundamental}. Thus, we can write $x' = j(z', a_5', a_{13}', \dots, a_{p_k}')$ where $a_5' \in 5^{\alpha'}\mathbb{Z}_5$ and
    \[|z'|_\mathbb{C} \le \frac{\eta}{2^4}, \quad |a_5'|_5 \le \frac{2^8}{\eta^2}\Xi_h, \quad \forall 2 \le j \le h, |a_{p_j}'|_{p_j} \le \frac{1}{p_j}, \quad \forall h < j \le k, |a_{p_j}'|_{p_j} \le 1.\]
    Let $a_5 \in 5^\ell\mathbb{Z}_5\setminus 5^{\ell+1}\mathbb{Z}_5$, so $\log_5(\eta^{-1}) \le \ell \le \log_5(k\varpi^{-1})$. Similarly, let $a_5' \in 5^{\ell'} \mathbb{Z}_5 \setminus 5^{\ell'+1}\mathbb{Z}_5$, so $\ell' \ge \alpha'$.

    First we wish to discard the case where $\ell' > \alpha'$. We achieve this as in the proof of Theorem 5.1 of \cite{kravitz2025}, by shifting the representative by some $j(\iota^\Delta(q))$ for $q \in A$ if necessary. Consider
    \[q := \overline{P_5}^{-\alpha'} \prod_{j = 2}^k \overline{P_{p_j}} \in A.\]

    Now we bound each absolute value individually. Note that $\alpha' \le \log_5(\eta^2\Xi_h^{-1}/2^8) + 1$. Hence
    \begin{align*}
        |z' - q|_\mathbb{C} &\le |z'|_\mathbb{C} + |q|_\mathbb{C} \\
        &\le \frac{\eta}{4} + (5^{\alpha'} \Xi_h)^{1/2} \\
        &\le \frac{\eta}{2^4} + \left(5 \cdot \frac{\eta^2}{2^8}\right)^{1/2} < \frac{\eta}{2}.
    \end{align*}
    Since $|a_5'|_5 < 5^{-\alpha'} = |q|_5$, we have from the ultrametric inequality that $|a_5' - q|_5 = 5^{-\alpha'}$. Furthermore, 
    \[\forall 2 \le j \le h, |a_{p_j}' - q'|_{p_j} \le \max\{|a_{p_j}'|_{p_j}, |q'|_{p_j}\} = \frac{1}{p_j}, \quad \forall h < j \le k, |a_{p_j}' - q|_{p_j} \le 1.\]
    Thus, after shifting our representative if necessary, we can assume $\ell' = \alpha'$. The cost this incurs is a weaker Archimedean bound $|z'|_\mathbb{C} \le \frac{\eta}{2}$.

    Now we push to the right scale, then cover every coprime residue modulo $T := \overline{P_5}^t$, where we choose $t := k^{1/80}$ with reference to our subset product coverage result, Theorem \ref{thm:subset-products}. In particular, there exists $0 \le s \le \ell - \alpha' \le \log_5(2^8k\varpi^{-1}\eta^{-2}\Xi_h)$ such that $|\theta_5^sa_5| = 5^{-\alpha'}$, and a further $m \in M_k$ such that
    \[|a_5' - \theta_5^sma_5|_5 \le 5^{-t-\alpha'} \le 5^{-\alpha}.\]
    We also have in the remaining non-Archimedean components,
    \[\forall 2 \le j \le h, |a_{p_j}' - \theta_5^sma_{p_j}|_{p_j} \le \max\{|a_{p_j}'|_{p_j}, p_j \cdot |a_{p_j}|_{p_j}\} \le \frac{1}{p_j},\]
    \[\forall h < j \le k, |a_{p_j}' - \theta_5^sma_{p_j}|_{p_j} \le \max\{|a_{p_j}'|_{p_j}, p_j \cdot |a_{p_j}|_{p_j}\} \le 1.\]
    Finally, in the complex component,
    \[|z' - \theta_5^smz|_{\mathbb{C}} \le |z'|_\mathbb{C} + |z|_{\mathbb{C}} \le \eta.\]
    This suffices once we notice that the complex component of $F(n, h, k, \alpha)$, a complex square of side length $(5^{\alpha - n}\Xi_h)^{1/2}$, contains the complex $\eta$-ball given the constraint on $\eta$. 
\end{proof}

One can repeat the same proof assuming a different scaling prime $p_j$. We note that the variation in the constraints is quantitatively insignificant in Section \ref{subsec:minor}. 

Note that we must have $h$ polynomially larger than $k$. Thus, our first set of mixing primes is significantly larger than the second. Indeed, without loss of generality, we obtain $|a_5'|_5 = 5^{-\alpha'} \le \Xi_h^{-1}$ and require $|a_5'|_5^{-1}$ to be a scale at which we can generate all the coprime residues using $M_k$. In particular, $\Xi_h \le \exp(k^{1/80})$.  

\subsubsection{Entropy input.} \label{subsec:entropy} Here we follow the blueprint set out in Section 5.2 of \cite{kravitz2025} and Section 3 of \cite{bourgain2009}. We recall a few pieces of notation, mostly taken from \cite{kravitz2025}. 

For $\mathcal{P}$ a finite partition of $X$ and $\mu$ a probability measure on $X$, we consider the \emph{entropy} of $\mathcal{P}$ with respect to $\mu$ as
\[H_\mu(\mathcal{P}) := -\sum_{D \in \mathcal{P}} \mu(D) \log(\mu(D))\]
with the convention that $0\log(0) = 0$. For $a \in A$, we denote by $a_*\mu$ the pushforward of $\mu$ under the $\times a$ map. Note that, for example, this gives
\[(a_*\mu)(f) = \int_X f(ax)d\mu(x).\]

Let $\tilde{\mu}$ denote the Haar measure on $X$. On the fundamental domain $F_0$, this is given by the product of the Haar measures on the individual components. Let $\lambda$ denote the Lebesgue measure on $[-1/2, 1/2)$. 

Let $f_\eps$ denote the smooth bump function adapted to the dual pyjama stripe $E^*(\eps)$ within the fundamental domain $F_0$, given by
\[f_\eps(x, y, a_5, a_{13}, \dots, a_{p_k}) = G_\eps(x) = g(x/\eps).\]
Here $g : \mathbb{R} \to \mathbb{R}_{\ge 0}$ is a smooth bump function satisfying
\[g(x) = \begin{cases}
    1 & |x| \le 1/10 \\
    0 & |x| \ge 1/5.
\end{cases}\]
Note that $\tilde{\mu}(f_\eps) = \lambda(G_\eps)$. We show that $f_\eps$ is actually a smooth bump function adapted to $E^*(\eps)$ in any complex translate of $F_0$. 

\begin{lmm}\label{lmm:f-translate}
    Let $F'$ denote an arbitrary complex translate of $F_0$, i.e. for $z \in \mathbb{C}$, 
    \[F' = (z + [-1/2, 1/2) \times [-1/2, 1/2)i) \times \mathbb{Z}_5 \times \dots \times \mathbb{Z}_{p_k}.\]

    Then $f_\eps$ is a smooth bump function adapted to $E^*(\eps)$ within $F'$. 
\end{lmm}
\begin{proof}
    Recall our argument for $E^*(\eps)$ being a prism-like set in $F_0$, at the start of Section \ref{subsec:fundamental}. Let
    \[\chi = j(z, a_5, \dots, a_{p_k}) = j(z', a_5', \dots, a_{p_k}')\]
    denote its representatives in $F_0, F'$ respectively. Then $-\Re(z) = \chi(1) = -\Re(z')$ and thus $f_\eps$ takes the exact same shape on $F'$ as on $F_0$, since it is a function of real part only. 
\end{proof}

We introduce Fourier analysis on $X$. As usual, let $e(t)= \exp(2\pi it)$, and let $w \in A$. We factorise $\widehat{f_\eps}$ in the fundamental domain $F_0$ as
\[\widehat{f_\eps}(w) = \left(\int_{-1/2}^{1/2} e(\Re(w) x) G_\eps(x)dx\right) \left(\int_{-1/2}^{1/2} e(-\Im(w) y) dy\right) \prod_{j = 1}^k \int_{\mathbb{Z}_{p_j}} e(-\{a_{p_j}w\}_{p_j})da_{p_j}.\]
Given the shape of our $f_\eps$, we expect our Fourier analysis to be virtually identical to Fourier analysis on $\mathbb{R}/\mathbb{Z}$. 

\begin{lmm}\label{lmm:fourier-vanish}
    Define $f_\eps : X \to \mathbb{R}$ as above. Then
    \[\widehat{f_\eps}(w) = \begin{cases}
        \widehat{G_\eps}(w) & w \in \mathbb{Z} \\
        0 & \text{otherwise.}
    \end{cases}\]
\end{lmm}
\begin{proof}
    Suppose $\widehat{f_\eps}(w) \neq 0$. Note that $\mathbb{Z}_{p_j}$ is invariant under translation, say by 1. Hence
    \begin{align*}
        \int_{\mathbb{Z}_{p_j}} e(-\{a_{p_j}w\}_{p_j})da_{p_j} &= \int_{\mathbb{Z}_{p_j}} e(-\{(a_{p_j} - 1)w\}_{p_j})da_{p_j} \\
        &= e(\{w\}_{p_j}) \int_{\mathbb{Z}_{p_j}} e(-\{a_{p_j}w\}_{p_j})da_{p_j}.
    \end{align*}
    Thus, if $\{w\}_{p_j} \not \in \mathbb{Z}$, the integral vanishes. Hence, since $\{w\}_{p_j} \in [0,1)$, we have $w \in \mathbb{Z}_{p_j}$ for all $1 \le j \le k$. We can then suppose $w \in \mathbb{Z}[i]$. 

    Similarly, we have that if $\Im(w) \neq 0$, then
    \[\int_{-1/2}^{1/2} e(-\Im(w)y)dy = 0.\]
    Hence, for $\widehat{f_\eps}(w) \neq 0$, we must have $w \in \mathbb{Z}$. In this case, we end up with
    \[\widehat{f_\eps}(w) = \int_{-1/2}^{1/2} e(w x) G_\eps(x) dx = \widehat{G_\eps}(w).\qedhere\]
\end{proof}

Since we work with Fourier analysis on $\mathbb{R}/\mathbb{Z}$, we can heavily simplify the arguments in Section 5.2 of \cite{kravitz2025}. In particular, we can almost work within the original framework of Bourgain, Lindenstrauss, Michel, Venkatesh \cite{bourgain2009}, using the standard $\ell^2$ norm instead of a Sobolev norm. 

First, we wish to obtain a result for discretised measures and functions on
\[\sigma + j(\iota^\Delta((\theta_5^{-\ell}A)/A)) \subseteq X,\]
which we abuse notation and simplify to $\sigma + (\theta_5^{-\ell}A)/A$. Let $\Omega_\ell := (\mathbb{Z}[i]/(P_5^\ell A))^\times \subseteq A/(\theta_5^\ell A)$. 

We consider a variation of both Lemma 5.4 of \cite{kravitz2025} and Lemma 3.5 of \cite{bourgain2009}, adopting the same proof structure. It states that on average over $\Omega_\ell$, $(\xi_*\mu)(f_\eps)$ is close to $\tilde{\mu}(f_\eps) = \lambda(G_\eps)$. 

\begin{lmm}\label{lmm:ell_2}
    Let $\ell$ be a positive integer. Let $\sigma \in X$ be arbitrary and $\mu$ a probability measure on $\sigma + (\theta_5^{-\ell}A)/A$. Then we have, for $f_\eps, G_\eps$ defined as above,
    \[\ev_{\xi \in \Omega_\ell} |(\xi_*\mu)(f_\eps) - \lambda(G_\eps)| \le C_5||G_\eps'||_2^2 \cdot ||\mu||_2^2,\]
    where $C_5 := 5 \cdot 2\frac{\pi^2}{3}$.
\end{lmm}
\begin{proof}
    First, we demonstrate that all of the non-trivial Fourier coefficients of $(\xi_*\mu)$ are small when averaging over $\Omega_\ell$. We have from Lemma \ref{lmm:dual} that the dual of $(\theta_5^{-\ell}A)/A$ is isomorphic to $A/(\theta_5^\ell A)$ where, given $w \in A/(\theta_5^\ell A)$, we have the character
    \[w(x) := -\Re(wx) + \sum_{j = 1}^k \{\iota_{p_j}(w \cdot x)\}_{p_j} \pmod{1}\]
    for $x \in (\theta_5^{-\ell}A)/A$. 
    
    Equipped with Lemma \ref{lmm:intertwine}, we can proceed exactly as in the first step of the proof of Lemma 5.4 of \cite{kravitz2025}. Furthermore, we can set (within the terminology of the proof) $\alpha = 1$. Indeed, $\Omega_\ell$ is the set of all coprime cosets, and so it is invariant under adding $\theta_5$ to itself. This yields
    \[\ev_{\xi \in \Omega_\ell} |\widehat{\xi_*\mu}(w)|^2 \le 5^{\min\{1 + v_{P_5}(w), \ell\}} ||\mu||_2^2.\]

    Next, we proceed as in the second step of the proof of Lemma 3.5 of \cite{bourgain2009}. We expand $f_\eps$ as a Fourier series, keeping only the integer modes by Lemma \ref{lmm:fourier-vanish},
    \[f_\eps(x) = \sum_{w \in \mathbb{Z}} \widehat{G_\eps}(w)e(w(x)),\]
    with $\widehat{f_\eps}(0) = \lambda(G_\eps)$. By applying Parseval's Identity and the Cauchy-Schwarz inequality,
    \begin{align*}
        \ev_{\xi \in \Omega} |(\xi_*\mu)(f_\eps) - \lambda(G_\eps)|^2 &= \ev_{\xi \in \Omega} \left|\sum_{0 \neq w \in \mathbb{Z}} \widehat{G_\eps}(w) \widehat{\xi_*\mu}(w)\right|^2 \\
        &\le \left(\ev_{\xi \in \Omega} \sum_{0 \neq w \in \mathbb{Z}} \frac{|\widehat{\xi_*\mu}(w)|^2}{w^2}\right) \left(\sum_w w^2 |\widehat{G_\eps}(w)|^2\right) \\
        &= \left(\sum_{0 \neq w \in \mathbb{Z}} w^{-2} \ev_{\xi \in \Omega} |\widehat{\xi_*\mu}(w)|^2\right) ||G_\eps'||_2^2 \\
        &\le ||G_\eps'||_2^2 \cdot ||\mu||_2^2 \sum_{0 \neq w \in \mathbb{Z}} \frac{5^{\min\{1 + v_{P_5}(w), \ell\}}}{w^2}.
    \end{align*}
    In the final step, we used our above bound on $\ev_{\xi \in \Omega_\ell} |\widehat{\xi_*\mu}(w)|^2$. 
    
    We can now bound the constants as follows. Note that since $w$ is an integer, $v_{P_5}(w) = v_5(w)$. 
    \begin{align*}
        \sum_{0 \neq w \in \mathbb{Z}} \frac{5^{\min\{1 + v_5(w), \ell\}}}{w^2} &\le \sum_{0 \le r \le \ell} \sum_{\substack{w \neq 0\\5^r \mid w}} \frac{5^{1+r}}{w^2} \\
        &\le \sum_{0 \le r \le \ell} \sum_{w \neq 0} \frac{5^{1+r}}{5^{2r}} \frac{1}{w^2} \\
        &\le \frac{\pi^2}{3} \sum_{0 \le r \le \ell} 5^{1-r} \le 5 \cdot 2\frac{\pi^2}{3}.\qedhere
    \end{align*}
\end{proof}

One can alternatively take an approach more close to Lemma 5.4 of \cite{kravitz2025} by using a suitable Sobolev norm to make the constant $C$ absolute (independent of choice of prime). Since the resultant loss from $C_{p_j} = p_j \cdot 2\frac{\pi^2}{3}$ ends up being (quantitatively) insignificant, we take this simpler approach.

Now we make use of our entropy assumption. If $\mathcal{P}(n, h, k, \alpha)$ has large entropy with respect to $\mu$, then some $(\theta_5^s)_*\mu$ majorises another measure $\nu$ that is quite spread out in an $\ell^2$ sense. This is exactly Lemma 5.5 of \cite{kravitz2025}, which follows from Lemma 3.6 of \cite{bourgain2009} since we have identical intertwining behaviour between $(\theta_5^{-n}A)/A \cong \mathbb{Z}/5^n\mathbb{Z}$. 

\begin{lmm}\label{lmm:spread-out}
    Let $\rho, \delta > 0$ be such that $\delta \le \rho/10$. Let $n$ be a positive integer. Suppose $\mu$ is a measure on $(\theta_5^{-n}A)/A$ with entropy at least $\rho \log N$ with respect to the partition into singletons. Then for each integer $10/\delta \le \ell \le \delta n$, we have $(\theta_5^s)_*\mu \ge \nu$ for some integer $0 \le s \le (1-\delta)n$ and measure $\nu = \sum_{i =1}^{5^{n-\ell}} w_i\nu_i$, such that
    \begin{outline}[enumerate]
        \1 the $\nu_i$'s are probability measures, each supported on a translate of $(\theta_5^{-\ell}A)/A$,
        \1 we have $w_i \in \mathbb{R}_{\ge 0}$ with $\sum_{i=1}^{5^{n-\ell}} w_i \ge \rho - 2\delta$,
        \1 overall we have $\sum_{i = 1}^{5^{n-\ell}} w_i ||\nu_i||_2^2 \le 2 \cdot 5^{-\delta \ell}$. 
    \end{outline}
\end{lmm}

This allows us to prove our effective result for our discretised measure. This follows a similar structure to Proposition 5.6 of \cite{kravitz2025} and Proposition 3.8 of \cite{bourgain2009}, but we chose our parameters to ensure quasipolynomial loss instead of logarithmic.

\begin{prop}\label{prop:discretised}
    Let $\rho \in (0.1,0.9)$ and $\delta := \rho/10$ be parameters. Let $n, h$ be sufficiently large positive integers with
    \[N := 5^n, \quad h := n^{1/2}.\]
    Suppose $\mu$ is a measure on $(\theta_5^{-n}A)/A$ with entropy at least $\rho \log N$ with respect to the partition into singletons. Then there exists an integer $0 \le s \le (1-\delta)n$ and some $m \in M_h$ such that, for $f_\eps, G_\eps$ defined above,
    \[(\theta_5^sm)_*\mu(f_\eps) \ge (\rho - 2\delta) \lambda(G_\eps) - 5^{-\delta h^{1/80}/4} ||G_\eps'||_2.\]
\end{prop}
\begin{proof}
    Set $\ell = n^{1/160}$. Hence, for sufficiently large $n$, we obtain $10/\delta \le \ell \le \delta n$. Applying Lemma \ref{lmm:spread-out}, we obtain a measure $\nu = \sum_{i = 1}^{5^{n-\ell}} w_i\nu_i$ alongside an integer $0 \le s \le (1-\delta)n$, such that $(\theta_5^s)_*\mu \ge \nu$, and $w := \sum_{i = 1}^{5^{n-\ell}} w_i \ge \rho - 2\delta$, $\sum_{i = 1}^{5^{n-\ell}} w_i||\nu_i||_2^2 \le 2 \cdot 5^{-\delta \ell}$.

    We can compute
    \begin{align*}
        \ev_{\xi \in \Omega_\ell} |(\xi_*\nu)(f_\eps) - w\lambda(G_\eps)| &\le \sum_{i=1}^{5^{n-\ell}} w_i \ev_{\xi \in \Omega_\ell} |(\xi_*\nu_i)(f_\eps) - \lambda(G_\eps)| \\
        &\le w^{1/2} \left(\sum_{i=1}^{5^{n-\ell}} w_i \left(\ev_{\xi \in \Omega_\ell} |(\xi_*\nu_i)(f_\eps) - \lambda(G_\eps)|\right)^2\right)^{1/2} \\
        &\le \left(\sum_{i=1}^{5^{n-\ell}} w_i \ev_{\xi \in \Omega_\ell} |(\xi_*\nu_i)(f_\eps) - \lambda(G_\eps)|^2\right)^{1/2} \\
        &\le \left(\sum_{i=1}^{5^{n-\ell}} w_i C_5 ||G_\eps'||_2^2 \cdot ||\nu_i||_2^2\right)^{1/2} \\
        &\le C_5^{1/2} ||G_\eps'||_2 \left(2 \cdot 5^{-\delta \ell}\right)^{1/2}.
    \end{align*}
    Note that the inequalities in order follow from: the triangle inequality after decomposing $\nu = \sum_{i = 1}^{5^{n-\ell}} w_i \nu_i$; the Cauchy-Schwarz inequality; $w \le 1$ and the Cauchy-Schwarz inequality; Lemma \ref{lmm:ell_2}; Lemma \ref{lmm:spread-out}, respectively. 

    Thus, there exists some $\xi \in \Omega_\ell$ such that
    \[|(\xi_*\nu)(f_\eps) - w\lambda(G_\eps')| \le (2C_5)^{1/2} ||G_\eps'||_2 \cdot 5^{-\delta \ell/2}.\]
    Using Theorem \ref{thm:subset-products} on subset products, we find that $M_h$ generates all coprime residues modulo $5^\ell$ for $\ell = h^{1/80}$. Thus, we can find $m \in M_h$ such that $m - \xi \in \theta_5^\ell A$, which gives
    \[|(m_*\nu)(f_\eps) - w\lambda(G_\eps)| \le (2C_5)^{1/2} ||G_\eps||_2 \cdot 5^{-\delta \ell/2}.\]
    Since $(\theta_5^s)_*\mu \ge \nu$, we find that
    \begin{align*}
        (\theta_5^sm)_*\mu(f_\eps) &\ge (\rho - 2\delta) \lambda(G_\eps) - (2C_5)^{1/2} ||G_\eps'||_2 \cdot 5^{-\delta \ell/2} \\
        &\ge (\rho - 2\delta)\lambda(G_\eps) - \kappa_5 5^{-\delta h^{1/80}/2} ||G_\eps'||_2
    \end{align*}
    where $\kappa_5 = (2C_5)^{1/2} = 5^{1/2} \cdot \frac{2\pi}{\sqrt{3}}$. Yet since $h$ is sufficiently large, we can (very loosely) bound
    \[\kappa_5 5^{-\delta h^{1/80}/2} \le 5^{-\delta h^{1/80}/4}. \qedhere\]
\end{proof}

One of our main new ideas reveals itself in this final step, converting from the discretised measure to our original one. Our specific choice of partition means that even after applying $\theta_5^sm$, we land in some complex translate of $F_0$. Lemma \ref{lmm:f-translate} tells us $f_\eps$ is determined by real part only. Hence, to bound the difference between these measures (using the triangle inequality), we only need to consider the real difference and ignore the $p$-adic components. 

\begin{thm}\label{thm:effective}
    Let $\rho \in (0.1,0.9)$ and $\delta := \rho/10$ be parameters. Let $n, h, k, \alpha$ be sufficiently large positive integers with
    \[2 \le h \le k, \quad N := 5^n, \quad h := n^{1/2}, \quad \alpha \ge (1-\delta)n.\]
    Suppose $\mu$ is a probability measure on $X$ such that
    \[H_\mu(\mathcal{P}(n, h, k, \alpha)) \ge \rho \log N.\]
    Then there is an integer $0 \le s \le (1-\delta) n$ and some $m \in M_h$ such that, for $f_\eps, G_\eps$ defined above,
    \[(\theta_5^sm)_*\mu(f_\eps) \ge (\rho - 2\delta) \lambda(f_\eps) - 5^{-\delta h^{1/80}/4} ||G_\eps'||_2 - (2 \cdot 5^{\alpha - n}\Xi_h)^{1/2} ||G_\eps'||_\infty.\]
\end{thm}
\begin{proof}
    We have that the parts of the partition $\mathcal{P}(n, h, k, \alpha)$ are indexed by the elements of $(\theta_5^{-n}A)/A$. Thus, we can define a probability measure $\mu'$ on $(\theta_5^{-n}A)/A$ given by
    \[\mu'(\{a\}) := \mu(D_a)\]
    where $D_a$ is the part of $\mathcal{P}(n, h, k, \alpha)$ containing $a$, a translation of $F(n, h, k, \alpha)$. By construction, the entropy of $\mu'$, with respect to the partition into singletons,  is precisely $H_\mu(\mathcal{P}(n, h, k, \alpha))$. Applying Proposition \ref{prop:discretised} to $\mu'$, we obtain integers $0 \le s \le (1-\delta) n$ and $m \in M_h$ such that
    \[(\theta_5^sm)_*\mu'(f_\eps) \ge (\rho - 2\delta) \lambda(G_\eps) - \kappa_55^{-\delta h^{1/80}/2} ||G_\eps'||_2,\]
    where $\kappa_5 > 0$ is the constant given at in the statement of Proposition \ref{prop:discretised}. 

    Now note that
    \[|(\theta_5^sm)_*\mu(f_\eps) - (\theta_5^sm)_*\mu'(f_\eps)| \le \sup_{\substack{\text{$x, x'$ same part}\\\text{of $\mathcal{P}(n,h,k,\alpha)$}}} |f_\eps(\theta_5^smx) - f_\eps(\theta_5^smx')|.\]

    Suppose $x, x' \in D_a \subseteq F_{h, k, \alpha}$, as noted at the end of Section \ref{subsec:fundamental}. We take representatives within $F_{h, k, \alpha}$, say
    \[x = j(z, a_5, \dots, a_{p_k}), \quad x' = j(z', a_5', \dots, a_{p_k}').\]
    Our choice of fundamental domains ensures
    \[a_5, a_5' \in 5^\alpha\mathbb{Z}_5, \quad \forall 2\le j\le h, a_{p_j}, a_{p_j}' \in p_j\mathbb{Z}_{p_j}, \quad \forall h < j \le k, a_{p_j}, a_{p_j}' \in \mathbb{Z}_{p_j}.\]
    Furthermore, since $x, x' \in D_a$, we have $x - x' \in F(n, h, k, \alpha)$ hence $|z - z'|_\mathbb{C} \le (2 \cdot 5^{\alpha - n}\Xi_h)^{1/2}$. Now consider
    \[\theta_5^smx = j(\theta_5^smz, \theta_5^sma_5, \dots, \theta_5^sma_{p_k}), \quad \theta_5^smx' = j(\theta_5^smz', \theta_5^sma_5', \dots, \theta_5^sma_{p_k}').\]
    Since $\alpha \ge (1-\delta)n \ge s$ and $m \in M_h$, we obtain
    \[\forall 1 \le j \le k, \theta_5^sma_{p_j}, \theta_5^sma_{p_j}' \in \mathbb{Z}_{p_j}.\]
    Thus, $\theta_5^smx, \theta_5^smx'$ lie in some complex translation translate of $F_0$. By Lemma \ref{lmm:f-translate}, $f_\eps$ only cares about real part. Hence, by the mean value theorem,
    \begin{align*}
        |f_\eps(\theta_5^smx) - f_\eps(\theta_5^smx')| &= |\Re(\theta_5^smz) - \Re(\theta_5^smz')| \cdot \frac{|G_\eps(\Re(\theta_5^smz)) - G_\eps(\Re(\theta_5^smz))|}{|\Re(\theta_5^smz) - \Re(\theta_5^smz')|} \\
        &\le |\theta_5^sm(z-z')|_\mathbb{C} \cdot ||G_\eps'||_\infty \\
        &\le (2 \cdot 5^{\alpha - n}\Xi_h)^{1/2} ||G_\eps'||_\infty.
    \end{align*}
    The triangle inequality then tells us
    \[(\theta_5^sm)_*\mu(f_\eps) \ge (\rho - 2\delta) \lambda(G_\eps) - \kappa_55^{-\delta h^{1/80}/2} ||G_\eps'||_2 - (2 \cdot 5^{\alpha - n}\Xi_h)^{1/2} ||G_\eps'||_\infty.\qedhere\]
\end{proof}

Now we pick $\eps > 0$ to obtain a result of the following form. Given a set which intersects at least $N^\rho$ parts of $\mathcal{P}(n, h, k, \alpha)$, we can run a short orbit to obtain quantitative denseness at a quasipolynomial scale. 

\begin{thm}\label{thm:blmv-final}
    Let $\rho \in (0.1,0.9)$ and $\delta := \rho/10$ be parameters.  Let $n, h, k, \alpha$ be sufficiently large positive integers with
    \[2 \le h \le k, \quad N := 5^n, \quad h := n^{1/2}, \quad (1-\delta)n \le \alpha \le (1-\delta/4)n.\]
    Suppose $\Gamma \subseteq X$ intersects at least $N^\rho$ parts of $\mathcal{P}(n, h, k, \alpha)$. Then
    \[\{\theta_5^smx: m \in M_h, 0 \le s \le S_2, x \in \Gamma\}, \quad S_2 = (1-\delta)n\]
    intersects $E^*(\eps_2)$ for $\eps_2 := 5^{-\delta n^{1/160}/20}$.
\end{thm}
\begin{proof}
    There exist fixed constants $c_1, c_2 > 0$ depending on $g$ only such that $\lambda(G_\eps) \ge c_1\eps$ and $||G_\eps'||_\infty \le c_2\eps^{-1}$, uniformly across $\eps$. Furthermore, $\rho - 2\delta > 0.08$. Thus, 
    \[(\theta_5^sm)_*\mu(f_\eps) \ge 0.08c_1\eps - c_2(5^{-\delta h^{1/80}/4} + (2 \cdot 5^{\alpha - n}\Xi_h)^{1/2}) \eps^{-1}.\]
    For an elementary bound, we have $\Xi_h \le p_h^h \le \exp(h^{3/2})$. Thus, we observe
    \begin{align*}
        5^{-\delta h^{1/80}/4} + (2 \cdot 5^{\alpha - n}\Xi_h)^{1/2} &\le 5^{-\delta n^{1/160}/4} + (2 \cdot 5^{-\delta n/4} \exp(n^{3/4}))^{1/2} \\
        &\le 5^{-\delta n^{1/160}/8}
    \end{align*}
    since $n$ is sufficiently large. In particular, pick $\eps = \eps_2$ where $\eps_2 := 5^{-\delta n^{1/160}/20}$. Then
    \[(\theta_5^sm)_*\mu(f_{\eps_2}) \ge 0.08c_1 5^{-\delta n^{1/160}/20} - c_2 5^{-3\delta n^{1/160}/20} > 0. \qedhere\]
\end{proof}

\subsubsection{Choosing parameters.} \label{subsec:minor} We collect all our parameters and constraints so far, and generalise from working with scaling prime 5 to an arbitrary $p_j$. We express all quantities in terms of $n$ for convenience. 

From Theorem \ref{thm:manners-denseness}, we have, noting that $\Xi_h$ is the largest possible product of primes that we need to take,
\[\alpha \le k^{1/80} + \log_{p_j}(\eta^2\Xi_h^{-1}/2^8), \quad \eta \le p_j^{(\alpha - n)/2} \Xi_h^{1/2} p_k^{-1/2}/2.\]
From Theorem \ref{thm:blmv-final}, we have
\[h := n^{1/2}, \quad (1-\delta)n \le \alpha \le (1-\delta/4)n.\]
Observe the elementary bounds $2^h \le \Xi_h \le \exp(h^{3/2})$. Thus, we can take for $n$ sufficiently large, given our constant $\delta \in (0.01, 0.09)$,
\[k := ((1+\delta/2)n)^{80}, \quad \eta := p_j^{-\delta n/2}, \quad h := n^{1/2}, \quad \alpha := (1-\delta/2)n -  O(n^{1/2}).\]

Recall from Theorem \ref{thm:manners-denseness}, we had our rotations in the Manners denseness step given by the product of
\[|M_k| = 2^k = 2^{((1 + \delta/2)n)^{80}}, \quad S_1 := \log_{p_j}(k\varpi^{-1}\eta^{-2} \Xi_h/2^8) = \log_{p_j} \varpi^{-1} + \delta n + O(n^{1/2}).\]
In the entropy denseness step in Theorem \ref{thm:blmv-final}, our rotations were the product of
\[|M_h| = 2^h = 2^{n^{1/2}}, \quad S_2 := (1-\delta)n.\]
This altogether intersects the dual pyjama stripe of width
\[\eps_2 := p_j^{-\delta n^{1/160}/20}.\]

To account for all such primes, consider for $T$ a positive integer,
\[\Theta_T := \left\{\prod_{j = 1}^k \theta_{p_j}^{e_j} : 0 \le e_1, \dots, e_k \le T, \text{at most one of $e_j$ is $> 1$}\right\}.\]
Note that $|\Theta_T| \le k2^kT$. We combine our two mixing steps in the next theorem, a simpler variation on Proposition 5.7 of \cite{kravitz2025}. For brevity, we define
\[\theta_5^{[0,S]} := \{\theta_5^s : 0 \le s \le S\}.\]

\begin{thm}\label{thm:minor}
    Let $0 < \varpi < 1/10$ and let $\Theta_0$ be an initial set of actions. Let $n, k$ be sufficiently large positive integers with $k := (1.024n)^{80}$. Suppose $p \in X$ and $\phi_1, \phi_2 \in \Theta_0$ satisfy $(\phi_1 - \phi_2)p = j(z, a_5, \dots, a_{p_k})$ with
    \[|z|_\mathbb{C} \le \eta/2, \quad \varpi \le \sum_{j = 1}^k |a_{p_j}|_{p_j} \le \eta/2.\]
    where $\eta := p_k^{-0.024n}$. Then for
    \[S_1 := \log \varpi^{-1} + 0.048n + O(n^{1/2}), \quad S_2 := 0.952n, \quad \eps_2 := \exp(-0.0024n^{1/160}),\]
    the additional orbit given by
    \[\Theta_{S_2} \Theta_{S_1} \Theta_0 \cdot p\]
    intersects $E^*(\eps_2)$.
\end{thm}
\begin{proof}
    Set $\delta := 0.048$. Without loss of generality, suppose $|a_5|_5 = \max_{1 \le j \le n} |a_{p_j}|_{p_j}$. Indeed, our Theorems \ref{thm:manners-denseness} and \ref{thm:blmv-final} apply in general with the choices of parameters provided above. The overall shape of these is preserved if we take $p_j = 5$, noting that $\log p_k \le 100\log n$ is insignificant on these scales. 
    
    We have that $|a_5|_5 \ge \varpi/k$ by the pigeonhole principle. For sufficiently large $n$, we have for all $1 \le j \le k$,
    \[\frac{1}{p_j^2} \ge \frac{1}{p_k^2} \ge k^{-3} \ge n^{-300} \ge p_j^{-\delta n/2}.\]
    Thus, for $h := n^{1/2}$, $\alpha := 0.976n - O(n^{1/2})$, and $\eta, S_1$ as given in the statement, we can apply Theorem \ref{subsec:manners} to $(\phi_1 - \phi_2)p$. By construction, the resultant orbit lies within
    \[\theta_5^{[0,S_1]}M_h \cdot (\phi_1 - \phi_2)p \subseteq \Theta_{S_1} \cdot (\Theta_0 - \Theta_0) p\]
    and hence this latter set intersects every set of the form $x' + F(n, h, k, \alpha)$. In particular, it intersects every part of the partition $P(n, h, k, \alpha)$. Our choice of superset $\Theta_{S_1}$ accounts for any arbitrary mixing prime $p_j$ in place of 5. 

    At the end of Section \ref{subsec:fundamental}, we showed that the difference of any two parts of $\mathcal{P}(n, h, k, \alpha)$ intersects four parts of $\mathcal{P}(n, h, k, \alpha)$. Since there are $N$ parts in total, we find that
    \[\Theta_{S_1} \Theta_0 \cdot p\]
    intersects $\gg N^{1/2} \ge N^{0.48}$ parts of our partition. 

    Apply Theorem \ref{thm:blmv-final} with $\rho := 0.48$, $n, h, k, \alpha$ given as above, and $\Gamma = \Theta_{S_1}\Theta_0 \cdot p$. The output orbit lies within
    \[\theta_5^{[0,S_2]}M_h \cdot \Theta_{S_1} \Theta_0 p \subseteq \Theta_{S_2} \cdot \Theta_{S_1}\Theta_0p.\]
    Hence, this latter set intersects $E^*(\eps_2')$ for $\eps_2' := 5^{-\delta n^{1/160}/20}$. Our choice of superset and $\eps_2$ in the statement then covers all possible mixing primes $p_j$, noting that $\delta = 0.048$.
\end{proof}

\subsection{Pyjama conversion} \label{sec:conversion}

In the final section of the upper bound proof, we take our results from the dual space $X$ and convert them back to the original Pyjama Problem setting. We find that the uncovered set is contained within small complex balls around sparse ``almost-lattices": lattices with a sub-lattice removed. The irrational trick takes care of these remaining points. The structure of this section is adapted from Section 7 of \cite{kravitz2025}. 

\subsubsection{Rationality lemma.} \label{subsec:rationality} We combine the major arcs and minor arcs in a dichotomy, presented as a rationality lemma similar to Theorem 1.2 of \cite{kravitz2025}. Recall that the pigeonhole step produces a point that is close to 0. The dichotomy then differentiates between the point being extremely close to $\overline{B}_1^\mathbb{C}(0)$ (major arc), or merely somewhat close (minor arc).

\begin{thm}\label{thm:rationality}
    For $\eps, \delta > 0$ sufficiently small, there exists a set of multipliers $\Theta \subseteq \mathbb{Q}(i)$ of size at most $\exp((\log \eps^{-1})^{40000}) \cdot \log \delta^{-1}$ satisfying the following. 

    For each continuous character $p \in j_\mathbb{C}(\mathbb{C})$, either: $\Theta \cdot p$ intersects $E^*(\eps)$; or there is a torsion point $q \in X$ of order at most $\exp((\log \eps^{-1})^{20000})$ such that $\dist(p, \overline{B}_{1/2}^\mathbb{C}(q)) \le \delta$. 
\end{thm}
\begin{proof}
    We wish to apply Proposition \ref{prop:pigeon} such that there exists $r \in \Theta_0 - \Theta_0$ of height at most $H$ with $\dist(rp, 0) \le \eta/2$. In particular, we obtain the constraints
    \[\eta \ge 2(k+2)N_0^{-1/(k+2)}, \quad H \ge N_0^2.\]
    From the major arc argument in Proposition \ref{prop:major}, there exists a constant $C > 0$ such that we also take $K := (H\eps_1^{-1})^C$ and have the constraint
    \[\varpi := \min\left\{\frac{\delta}{(k+1)H^2}, \exp(-10K)\right\}.\]
    Hence, for sufficiently large $n$, we can set
    \[N_0 := \exp(n^{100}) \ge (p_k^{-0.024n}(k+2)^{-1}/2)^{-(k+2)}, \quad H := \exp(2n^{100}), \quad \eps_1 := \exp(-n^{100}),\]
    thus giving
    \[K = \exp(3Cn^{100}), \quad \varpi = \min\{\delta \exp(-5n^{100}), \exp(-10\exp(3Cn^{100}))\}.\]

    With these parameters established, we present the dichotomy. Let $p \in j_\mathbb{C}(\mathbb{C})$ denote a continuous character. 

    {\bf Initial pigeonhole}. By Proposition \ref{prop:pigeon}, we obtain an initial set of multipliers $\Theta_0$ of cardinality $N_0$ such that there is $r := \phi_1 - \phi_2 \in \Theta_0 - \Theta_0$ of height at most $H$, satisfying $\dist(rp, 0) \le \eta/2$. In particular, we can write $rp = j(z, a_5, \dots, a_{p_k})$ where
    \[|z|_\mathbb{C} + \sum_{j = 1}^k |a_{p_j}|_{p_j} \le \eta/2.\]

    {\bf Major arcs}. Suppose we have
    \[\sum_{j = 1}^k |a_{p_j}|_{p_j} \le \varpi.\]
    Then proposition \ref{prop:major} implies that either: $\dist(p, \overline{B}_{1/2}^\mathbb{C}(q)) \le \delta$ for some torsion point $q \in X$ of order at most $H^2$; or $\theta_5^{[0,K]}p$ intersects $E^*(\eps_1)$. Note that this orbit is contained in
    \[\Theta_Kp.\]

    {\bf Minor arcs}. Alternatively, suppose that
    \[\varpi \le \sum_{j = 1}^k |a_{p_j}|_{p_j} \le \eta/2.\]
    Then Theorem \ref{thm:minor} tells us that
    \[\Theta_{S_1} \Theta_{S_2} \Theta_0p\]
    intersects $E^*(\eps_2)$. 

    It now suffices to introduce $\eps$ into our parameters. For $n$ sufficiently large, we can take $\eps \ge \max\{\eps_1, \eps_2\} = \exp(-0.0024n^{1/160})$. Thus, certainly we can let
    \[n = (\log \eps^{-1})^{180}.\]

    Note that $S_1 \ge \log \varpi^{-1} \ge 10K$, and so $\Theta_K \subseteq \Theta_{S_1}$. Thus, we can take $\Theta := \Theta_{S_1}\Theta_{S_2}\Theta_0$ to be our set of multipliers, satisfying
    \begin{align*}
        |\Theta| = |\Theta_{S_1}| \cdot |\Theta_{S_2}| \cdot |\Theta_0| &\le k2^kS_1 \cdot k2^kS_2 \cdot N_0 \\
        &\le n \cdot 2^{2(1.024n)^{80}} \cdot (\log \varpi^{-1} + 0.048n + O(n^{1/2})) \cdot 0.952n \cdot \exp(n^{100}) \\
        &\le \exp(n^{120}) \cdot \max\{\log \delta^{-1} + 5n^{100}, 10\exp(3Cn^{100})\} \\
        &\le \exp(n^{120}) \cdot (\log \delta^{-1} + \exp(n^{120})) \le \exp(n^{200}) \cdot \log \delta^{-1}.
    \end{align*}
    Thus, we can conclude that, for sufficiently small $\eps > 0$,
    \[\text{torsion point order} \le H^2 \le \exp(2n^{100}) \le \exp((\log \eps^{-1})^{20000}),\]
    \[\text{number of multipliers} \le \exp((\log \eps^{-1})^{40000}) \cdot \log \delta^{-1}.\qedhere\]
\end{proof}

\subsubsection{Almost-lattice in the plane.} As in Section 7.1 of \cite{kravitz2025}, we consider the following neighbourhood of an almost-lattice. 
\[F(m, D) := \left(\frac{D}{m}\mathbb{Z}[i] \setminus D\mathbb{Z}[i]\right) + \overline{B}_{2/3}(0) \subseteq \mathbb{C}.\]

Let $M := \exp((\log \eps^{-1})^{20000})$. We consider the following variation of Lemma 7.1 of \cite{kravitz2025}. This converts the rationality lemma within $X$ to a rationality lemma in the plane, with low-order torsion points mapping to rational obstructions. 

\begin{thm}\label{thm:almost-lattice}
    There exists a constant $\tilde{C} > 0$ such that the following holds. For sufficiently small $\eps > 0$ and $R \ge \exp(\eps^{-\tilde{C}})$, there is a Gaussian integer $D$ with $R \le |D| \le 5R$ such that
    \[\bigcup_{1 \le m \le M} F(m, D) \cup \bigcup_{\theta \in \Theta'} E_\theta(\eps) = \mathbb{C},\]
    where we choose $\delta := \exp(-3\log R)$ for the set $\Theta$ given by Theorem \ref{thm:rationality}, and $\Theta' := \theta_5^{[0, \floor{\log_5R} + 1]}\Theta$.  
\end{thm}
\begin{proof}
    Apply Theorem \ref{thm:rationality} to obtain
    \[\bigcup_{\substack{q \in X\\\text{ord}(q) \le M}} \left(\overline{B}_{1/2}^\mathbb{C}(q) + \overline{B}_\delta(0)\right) \cup \bigcup_{\theta \in \Theta'} E_\theta^*(\eps) \supseteq j_\mathbb{C}(\mathbb{C}).\]
    Now we can apply $j_\mathbb{C}^{-1}$ to both sides and find that
    \[\bigcup_{\substack{q \in X\\\text{ord}(q) \le M}} j_\mathbb{C}^{-1}\left(\overline{B}_{1/2}^\mathbb{C}(q) + \overline{B}_\delta(0)\right) \cup \bigcup_{\theta \in \Theta'} E_\theta(\eps) = \mathbb{C}.\]
    We wish to show that
    \[\bigcup_{\substack{q \in X\\\text{ord}(q) \le M}} j_\mathbb{C}^{-1} \left(\overline{B}_{1/2}^\mathbb{C}(q) + \overline{B}_\delta(0)\right) \setminus \bigcup_{\theta \in \theta_5^{[0, \floor{\log_R} + 1]}} E_\theta(\eps) \subseteq \bigcup_{1 \le m \le M} F(m, D).\]

    Let $w$ be an element of the left-hand side. Then we can write
    \[w = j_\mathbb{C}^{-1}(q + x + y), \quad j_\mathbb{C}(w - x') = q + y,\]
    for $q$ a torsion point of order $m \le M$, $x = j_\mathbb{C}(x')$ with modulus at most 1/2, and $y \in X$ with $\dist(y, 0) \le \delta$. Multiplying by $m$, we get
    \[j_\mathbb{C}(m(w-x')) = mq + mz = mz.\]

    Set $D := \overline{P_5}^{\floor{\log_5R}+1}$, with $R \le |D| \le 5R$. We wish to detect proximity to $\frac{D}{m}\mathbb{Z}[i]$ by testing on
    \[\theta_5^{[0, \floor{\log_5R} + 1]} \cup i\theta_5^{[0, \floor{\log_5R} + 1]}.\]
    Indeed, for $\omega \in \{0,1\}$ and $0 \le s \le \floor{\log_5R} + 1$, we have
    \[|mz(i^\omega \theta_5^s)| = m|z(i^\omega\theta_5^s)| \le M \delta \cdot (1 + 5^{\floor{\log_5R}+1}) \le 1/10\]
    for sufficiently small $\eps > 0$, by the choice of $\delta$ (note $M$ is only quasipolynomial in $\eps^{-1}$). Thus, 
    \[||\Re(m(w-x')\theta_5^s)||_{\mathbb{R}/\mathbb{Z}} \le 1/10, \quad ||\Im(m(w-x')\theta_5^s)||_{\mathbb{R}/\mathbb{Z}} \le 1/10.\]
    In particular, $m(w-x')\theta_5^s$ lies within $\sqrt{2}/10$ of a Gaussian integer. 

    We have the following Claim from Section 5.3 of \cite{manners2013}, reappearing just before Lemma 7.1 of \cite{kravitz2025}. Let $0 < \eta < 1/\sqrt{5}$ and $T \in \mathbb{N}$: if $z \in \mathbb{C}$ is such that $\theta_5^sz$ lies within distance $\eta$ of a Gaussian integer for all $0 \le s \le T$, then $z$ itself lies within $\eta$ of $\overline{P_5}^T\mathbb{Z}[i]$. 

    Applying the Claim, we deduce that $m(w-x')$ lies within $\sqrt{2}/10$ of $D\mathbb{Z}[i]$. We can divide by $m$ to find that
    \[w - x' \in \frac{D}{m}\mathbb{Z}[i] + \overline{B}_{\sqrt{2}/10}(0).\]
    Since $|x'| \le 1/2$, we find that $w \in \frac{D}{m}\mathbb{Z}[i] + \overline{B}_{2/3}(0)$ since $1/2 + \sqrt{2}/10 < 2/3$. 

    Now it suffices to show $w \not \in D\mathbb{Z}[i] + \overline{B}_{2/3}(0)$. Suppose for the sake of contradiction $w = x + y$ for $x \in D\mathbb{Z}[i]$ and $|y| \le 2/3$. The quantitative approximation lemma, Lemma \ref{lmm:quantapprox}, tells us that there exists a constant $c > 0$ for which, for all $\eps > 0$,
    \[\{\theta_5^s : 0 \le s \le \eps^{-2/c}\}\]
    is $\eps^2$-dense in $S^1$. For sufficiently small $\eps > 0$, this tells us
    \[\{\Re(\theta_5^sy) : 0 \le s \le \eps^{-2/c}\}\]
    is $\eps/10$-dense in the interval $[-|y|, |y|]$, noting that $|y| \le 2/3$. 

    Set $\tilde{C} := 2/c$. Then for all $0 \le s \le \eps^{-2/c} \le \floor{\log_5R} + 1$, we have that $\theta_5^s = P_5^s/\overline{P_5}^s$ is a Gaussian rational with denominator dividing $D$. Hence, there exists $0 \le s \le \eps^{-2/c}$ with
    \[\Re(\theta_5^sw) = \Re(\theta_5^sx) + \Re(\theta_5^sy) = \Re(\theta_5^sy) \in (-\eps, \eps).\]
    Yet we have $w \not \in \bigcup_{\theta \in \theta_5^{[0, \floor{\log_5R} + 1]}} E_\theta(\eps)$, which is a contradiction.
\end{proof}

\subsubsection{Irrational trick.} \label{subsec:irrational} Now we consider the geometric trick of Zhu from Section 7.2 of \cite{kravitz2025}, generalising the irrational trick of Malikiosis, Matolcsi and Ruzsa \cite{malikiosis2013}. This adds some cleverly-chosen irrational rotations to our rotation set, which allows us to cover the remaining balls around the sparse almost-lattice. 

For $m \in \mathbb{N}$, set
\[\alpha_m = \frac{1}{2\pi}\arccos\left(\frac{1}{2m}\right),\]
satisfying $e(\alpha_m) + e(-\alpha_m) = 1/m$. For $M \in \mathbb{N}$, we consider additional rotations
\[V(M) := \{e(\pm \alpha_m) : 1 \le m \le M\}\cup \{1\}.\]
We have that $|V(M)| = 2M+1$. We directly quote Lemma 7.2 from \cite{kravitz2025}. 

\begin{lmm}\label{lmm:geometric}
    Let $M \in \mathbb{N}$. Then for every partition $V(M) = \bigsqcup_{m = 1}^M V_m$, there are integers $\lambda_v$, for $v \in V(M) \setminus \{1\}$, such that the following holds:
    \begin{outline}[enumerate]
        \1 $\sum_{v \in V \setminus\{1\}} \lambda_v v = 1$, 
        \1 for each $1 \le m \le M$, if $v \in V_m$, then $\lambda_v$ is a multiple of $m$,
        \1 $|\lambda_v| \le M^8$ for all $v \in V \setminus \{1\}$. 
    \end{outline}
\end{lmm}

From this, we deduce quasipolynomial pyjamas. 

\begin{proof}[Proof of Theorem \ref{thm:pyjamas}]
    Apply Theorem \ref{thm:almost-lattice} with $R = 5^M$. Set $V := V(M)$, then we take
    \[\tilde{\Theta} = V\theta_5^{[0, \floor{\log_5R} + 1]}\Theta\]
    to be our full set of rotations. Indeed, we claim that
    \[\bigcup_{\theta \in \tilde{\Theta}} E_\theta(\eps) = \mathbb{C}.\]
    Suppose for the sake of contradiction that
    \[z \in \mathbb{C} \setminus \bigcup_{\theta \in \tilde{\Theta}} E_\theta(\eps).\]
    For all $v \in V$, by Theorem \ref{thm:almost-lattice} applied to $vz$, we have $vz \in \bigcup_{1 \le m \le M} F(m, D)$. We can classify
    \[V_m := \{v \in V : vz \in F(m, D)\}.\]
    Certainly $\bigcup_{1 \le m \le M} V_m = V$, then we can arbitrarily trim elements in the $V_m$ until this is a disjoint union. 

    Apply Lemma \ref{lmm:geometric} to obtain integers $\lambda_v$ for $v \in V(M) \setminus \{1\}$. Multiplying the first property by $z$, we find that
    \[z = \sum_{v \in V\setminus \{1\}} \lambda_v vz.\]
    For $v \in V_m$, we find that $mvz \in D\mathbb{Z}[i] + \overline{B}_{2m/3}(0)$, say $mvz = J_v + E_v$ for $J_v \in D\mathbb{Z}[i]$ and $|E_v| \le 2m/3$. Then we find that
    \begin{align*}
        \left|z - \sum_{v \in V \setminus \{1\}} \frac{\lambda_v}{m} J_v\right| &= \left|\sum_{m = 1}^M \sum_{v \in V_m \setminus \{1\}} \left(\lambda_v vz - \frac{\lambda_v}{m}J_v\right)\right| \\
        &\le \sum_{m = 1}^M \sum_{v \in V_m\setminus \{1\}} \left|\frac{\lambda_v E_v}{m}\right| \\
        &\le \frac23 \sum_{v \in V \setminus \{1\}} |\lambda_v| \le \frac{4}{3} M^9,
    \end{align*}
    where we used the third property in the last inequality. We use the second property to notice that
    \[\sum_{v \in V\setminus \{1\}} \frac{\lambda_v}{m} J_v \in D\mathbb{Z}[i].\] 
    Thus, for $\eps > 0$ sufficiently small, $z$ is certainly within $\frac{5^M}{2M} \le \frac{|D|}{2M}$ of $D\mathbb{Z}[i]$. The closest point to 0 within $F(m, D)$ is distance at least $|D|/m - 2/3$. Hence, we find that
    \[z \not \in \bigcup_{1 \le m \le M} F(m, D),\]
    which yields a contradiction. 

    Finally, it suffices to bound $|\tilde{\Theta}|$. We have, for sufficiently small $\eps > 0$,
    \begin{align*}
        |\tilde{\Theta}| &\le |V| \cdot (\floor{\log_5R} + 2) \cdot |\Theta| \\
        &\le (2M+1) (M+2) \exp((\log \eps^{-1})^{40000}) (3\log R) \\
        &\le \exp((\log \eps^{-1})^{50000}). \qedhere
    \end{align*}
\end{proof}

\section{Lower bound} \label{lower}

In contrast to the upper bound, the ideas for the lower bound are elementary and live within the original Euclidean plane. Here, we treat the pyjama problem purely as one in combinatorial geometry. We restart our labelling of parameters afresh (other than $\eps$ itself). 

\subsection{Covering $S_{1-\eps}(0)$}\label{sec:first-circle}

We begin with an initial bound that forms the basis of our argument. 

\begin{prop}\label{prop:first-circle}
    Let $0 < \eps < 1/10$. If $\Theta$ is a finite set of unit-norm complex numbers such that $\bigcup_{\psi \in \Theta} E_\psi(\eps) = \mathbb{C}$, then we must have
    \[|\Theta| \ge \left(\frac{\pi}{2} + o(1)\right) \eps^{-1}.\]
\end{prop}
\begin{proof}
    Define $C_1 := S_{1-\eps}(0)$ for convenience. Let the \emph{central strip} be given by
    \[S := \{z \in \mathbb{C} : \Re(z) \in (-\eps, \eps)\}.\]
    Note that $C_1$ has circumference $2(1-\eps)\pi$. It also has the property that given any rotation $\theta$, the rotated pyjama stripe intersects $C_1$ exactly at two antipodal arcs within (the rotation of) $S$. In particular,
    \[E_\theta(\eps) \cap C_1 = (\theta^{-1}S) \cap C_1 = \theta^{-1} (\gamma^+ \cup \gamma^-),\]
    where $\gamma^+, \gamma^-$ are the topmost and bottommost arcs of $C_1$, centred on the $y$-axis, cut out by $S$. 

    To compute the length of $\gamma^+$, note that the angle $\alpha$ subtended by the arc satisfies $\sin(\alpha/2) = \eps/(1-\eps)$. Thus, letting $\ell$ denote length, we find that
    \[\ell(\gamma^+) = 2(1-\eps)\arcsin\left(\frac{\eps}{1-\eps}\right).\]

    One must cover $C_1$ in order to cover the plane. Hence, the number of rotations required is certainly at least
    \[\frac{2(1-\eps)\pi}{\ell(\gamma^+ \cup \gamma^-)} = \frac\pi2 \arcsin\left(\frac{\eps}{1-\eps}\right)^{-1} = \left(\frac{\pi}{2} + o(1)\right) \eps^{-1}.\qedhere\]
\end{proof}

\subsection{Covering a neighbourhood of $S_{2-\eps}(0)$} \label{sec:second-circle}

The motivation for the improvement is loosely as follows. Suppose we take a sets of rotations only slightly larger than the bound given in Proposition \ref{prop:first-circle}. Analysing the proof, we find that our rotations must be approximately well-spaced on the unit circle. We prove that such a configuration cannot cover a tiny neighbourhood of $S_{2-\eps}(0)$. 

We set out some initial notation. Define
\[x_0 := \frac{\pi}{2} \arcsin\left(\frac{\eps}{1-\eps}\right)^{-1}, \quad N_0 := \ceil{x_0}\]
to be the number of rotations needed to cover $S_{1-\eps}(0)$. Note that $x_0^{-1} = \frac{2}{\pi}\eps + O(\eps^2)$. We consider, alongside the central strip, the \emph{side strips} given by
\[S^\pm := \{z \in \mathbb{C} : \Re(z) \in \pm (1-\eps, 1+\eps)\}.\]

\subsubsection{Alternative perspective.} \label{subsec:perspective} Let $N := \ceil{N_0(1+\delta)}$ for $0 < \delta < \delta_0$. Suppose $\Theta$ is a set of rotations of cardinality at most $N$. Without loss of generality, we can add another arbitrary rotation if necessary to make $N$ even. We now consider our set of rotations as angles instead of unit modulus complex numbers, in particular we set
\[\Theta := \{e^{i\alpha_1}, \dots, e^{i\alpha_N}\}\]
where $\alpha_j \in [0, 2\pi)$ for $1 \le j \le N$. 

In this framework, rotations of the pyjama stripe $E(\eps)$ by $\Theta$ does not cover the plane if and only if there exists $z \in \mathbb{C}$ such that $e^{-i\alpha_j} z \not \in E(\eps)$ for all $1 \le j \le N$. In fact, suppose rotations of $E(\eps)$ by $\Theta$ cover $S_{1-\eps}(0)$. Then we find such a $z$ with
\[2 - (8\sqrt{3} + 1)\eps \le |z| \le 2 - \eps\]
which satisfies this property. Thus, since all the points we ever consider lie in $B_{2-\eps}(0)$, we may restrict our attention to
\[E(\eps) \cap B_{2-\eps}(0) = S^* \cap B_{2-\eps}(0), \quad S^* := S \cup S^+ \cup S^-.\]

We work in the perspective of rotating points away from the pyjama stripe. In particular, we let our \emph{points} be
\[R \cdot \{e^{-i\alpha_1}, \dots, e^{-i\alpha_N}\} \subseteq C := R \cdot S^1,\]
for a given $R = 2 - t\eps$, where $t \in (1, 8\sqrt{3} + 1)$. For some $R$ in this range, we find a \emph{revolution} $\omega \in S^1$ (to distinguish from our pyjama stripe rotations) of these points such that they all lie outside $S^*$, i.e.
\[R\omega \cdot \{e^{-i\alpha_1}, \dots, e^{-i\alpha_N}\} \cap S^* = \varnothing.\]

We let a \emph{gap} denote an arc between consecutive points. Note that this is different to the intersection of the individual strips with $C$. 

We can compute the length of each of the four (symmetric) arcs of $S^\pm \cap C$, denoted by \emph{side arcs}. These are $(2-t\eps)\beta$ where $\beta$ is the angle of the arc. By the mean value theorem (and trigonometry), we have
\begin{align*}
    \beta &= \arccos\left(\frac{1-\eps}{2-t\eps}\right) - \arccos\left(\frac{1+\eps}{2-t\eps}\right) \\
    &= \frac{2\eps}{2-t\eps} \frac{1}{\sqrt{1 - \left(\frac{1-\eps}{2-t\eps}\right)^2}} + O(\eps^2) \\
    &= \frac{2}{\sqrt{3}} \eps + O(\eps^2).
\end{align*}

Thus, the arcs themselves are length $\frac{4}{\sqrt{3}}\eps + O(\eps^2)$.

We make a crucial symmetry observation. Note that the pyjama stripes are invariant under reflection in the origin, i.e.
\[e^{-i\alpha}z \in E(\eps), S^* \text{ if and only if } e^{i\alpha}z \in E(\eps), S^*\]
respectively. 

Thus, reflecting points in the origin if necessary and removing duplicates, we can without loss of generality start with each of $e^{-i\alpha_j}$ on the right semicircle. Then the condition to cover $C_1$ is that the angle between consecutive points (other than the start and end) is at most $\pi/x_0$. We can then reflect every second point to enforce that the angle between all pairs of consecutive $\alpha_j$ (within $[0, 2\pi)$ with $0, 2\pi$ identified) is at most $2\pi/x_0 = 4\eps + O(\eps^2)$. 

\subsubsection{Average-looking configuration.} \label{subsec:avg-looking} Consider any $R = 2 - t\eps$ in the range, noting that $t$ is bounded. We wish to apply the pigeonhole principle to approximate an average-looking configuration. We cannot do this globally since $\Theta$ was arbitrary, but we can locally, in particular in the six neighbourhoods where $S^*$ intersects $C$. 

Note that $C$ has circumference $2(2 - t\eps)\pi$. Thus, the expected length of a gap is
\[\lambda := \mathbb{E}[\text{gap length}] = \frac{2(2-t\eps)\pi}{N} = \frac{8}{1+\delta} \eps + O(\eps^2).\]

We call a gap \emph{short} if its length is at most $(1-25\delta)\lambda$, and long otherwise. Note the property that if we shrink or dilate $R$, shortness and longness are preserved. 

Furthermore, call a gap \emph{special} if it intersects $S^*$ or is adjacent to one that does.

\begin{lmm}\label{lmm:avg-looking}
    Let $0 < \eps < 1/10$ and $N \in \mathbb{N}$ such that $N = \ceil{(1+\delta)N_0}$, where $0 < \delta < 0.01$ and $N_0 = \ceil{x_0}$ for
    \[x_0 = \frac{\pi}{2} \arcsin\left(\frac{\eps}{1-\eps}\right)^{-1}.\]
    Let $\Theta$ be a set of $N$ rotations such that the angle between consecutive rotations is at most $2\pi/x_0$, and let $C = S_{2-t\eps}(0)$ where $t$ is bounded. Then there exists a revolution of our points such that all special gaps are long. 
\end{lmm}
\begin{proof}
    The condition on consecutive rotations implies that each gap has length at most $2(2-t\eps)\pi/x_0$.
    
    Suppose a $\rho$ proportion of the $N$ gaps are short. Then:
    \begin{align*}
        \text{total length of gaps} &= 2(2-t\eps)\pi \\
        \Rightarrow \sum_{\text{short gap $\sigma$}} \ell(\sigma) + \sum_{\text{long gap $\sigma$}} \ell(\sigma) &= 2(2-t\eps)\pi \\
        \Rightarrow (\rho N) \cdot (1-25\delta)\lambda + ((1-\rho)N) \cdot \frac{2(2-t\eps)\pi}{x_0} &\ge 2(2-t\eps)\pi \\
        \Rightarrow \rho(1-25\delta) + (1-\rho)\frac{N}{x_0} &\ge 1.
    \end{align*}
    Since $N/x_0 = 1+\delta + O(\eps)$, we find that $\rho \le \frac{1}{26} + O(\eps)$. \\

    Consider a revolution chosen uniformly at random. As we revolve, each gap intersects $S^*$ in one of the two \emph{central arcs}, of length $2\eps + O(\eps^2)$, or four side arcs. Thus, (excessively) accounting for overhang, each short gap is special for a total proportion of the circumference of $C$ bounded above by
    \[m := \frac{1}{2(2-t\eps)\pi} \left(2 \cdot (2\eps + O(\eps^2)) + 4 \cdot \left(\frac{4}{\sqrt{3}}\eps + O(\eps^2)\right) + (12 \cdot 4) \cdot \frac{2(2-t\eps)\pi}{x_0}\right)\]
    \[= \frac{1}{\pi} \left(1 + \frac{4}{\sqrt{3}} + 48 \cdot \frac{1}{1+\delta}\right)\eps + O(\eps^2).\]
    Let this quantity be $m$. Thus, by linearity of expectation, using the fact that $\delta \le 0.01$ we get
    \begin{align*}
        \mathbb{E}[\#\text{special short gaps}] &\le \rho N \cdot m \\
        &= \frac{1}{26} \frac{N}{\pi} \eps \left(1 + \frac{4}{\sqrt{3}} + 48 \cdot \frac{1}{1+\delta}\right) + O(\eps^2) \\
        &= \frac{1}{26 \cdot 2} \left(1.01 \cdot \left(1 + \frac{4}{\sqrt{3}}\right) + 48\right) < 1.
    \end{align*}
    Thus, by averaging, there exists a revolution where all special gaps are long. 
\end{proof}

\subsubsection{Reduction to an interval problem.} \label{subsec:reduction} In an average-looking configuration, our central and side arcs are much smaller compared to the special gaps themselves. The hope is that our two degrees of freedom in revolution and radius are enough to exploit this. 

We work within the revolution given by Lemma \ref{lmm:avg-looking} where all special gaps are long. First, we focus on the central strip.

Consider the top central arc of $S \cap C$. If it crosses two long gaps, revolve until it lies fully within a long gap. Now the bottom central arc could possibly be split between two long gaps. Yet by the pigeonhole principle, one of these yields an interval of length at least $\frac12 (1-25\delta)\lambda$ for which both the top and bottom central arcs lie entirely within a long gap, i.e. do not contain a point. Call such an interval \emph{permissible}. \\

Our plan forwards is as follows. We interpret our two degrees of freedom as two \emph{operations}. The first revolves the points, keeping within the permissible interval of length
\[L \ge \frac12 (1 - 25\delta)\lambda = 4\frac{1-25\delta}{1+\delta}\eps + O(\eps^2).\]
We set $\tau := \frac{1-25\delta}{1+\delta}$. The second operation is to shrink $R$ between $2 - \eps$ and $2 - (8\sqrt{3} + 1)\eps$. 

We demonstrate that shrinking has negligible effect on points very close to $S$. Indeed, suppose $z \in C$ with $\Re(z) = \eps'$. Then if shrinking brings $z$ into $S$, we must certainly have
\[\frac{2 - (8\sqrt{3}+1)\eps}{2 - \eps} \eps' = \Re\left(\frac{2-(8\sqrt{3} + 1)\eps}{2-\eps}z\right) < \eps.\]
But this yields $\eps' = \eps + O(\eps^2)$. Hence, for some constant $K > 0$, if we ensure we are always $K\eps^2$ clear of the boundaries of the permissible interval, then the central strip can be safely ignored. \\

We first consider a simplification of our scenario to provide all the motivation without having to consider error terms. 

\begin{lmm}\label{lmm:reduction}
    Let $0 < \eps < 1/10$, $N \in \mathbb{N}$ with $N = (1+\delta)\frac{\pi}{2}\eps^{-1}$. Suppose the special gaps all have length $8\eps$, and side arcs exactly $\frac{4}{\sqrt{3}}\eps$. Furthermore, assume the following statement about intervals holds. 

    ``There exists a constant $\eta > 0$ such that the following holds. Let $x, y, z, w \in \mathbb{R}/\mathbb{Z}$. Let $I_x = (x - \frac{1}{4\sqrt{3}}, x + \frac{1}{4\sqrt{3}}) \pmod{1}$ and similarly define $I_y, I_z, I_w$. Then one can apply the following operations until $[0, \eta] \cap (I_x \cup I_y \cup I_y \cup I_z) = \varnothing$. 
    \begin{outline}[enumerate]
        \1 Increase $x, z$ and decrease $x, w$ uniformly by $t \in [0, 1/4]$.
        \1 Increase $x, y, z, w$ uniformly by $t \in [0,1)$."
    \end{outline}

    Then there exists a choice of revolution and radius such that none of our points lie in $E(\eps)$. 
\end{lmm}
\begin{proof}
    Note that the side arcs now take a $\frac{1}{2\sqrt{3}}$ proportion of the surrounding special gaps. Restrict to, without loss of generality, to the neighbourhood of the top left side arc. We wish to analyse the relative position of the side arcs with respect to their surrounding special gaps. 

    \begin{figure}[h!]
        \centering
        \begin{tikzpicture}
            \draw (0,0) circle (4cm);
            \draw (-2.5,4.5) -- (-2.5,-4.5);
            \draw (-1.5,4.5) -- (-1.5,-4.5);
            \draw (1.5,4.5) -- (1.5,-4.5);
            \draw (2.5,4.5) -- (2.5,-4.5);
            \draw[dotted] (2,4.5) -- (2,-4.5);
            \draw[dotted] (-2,4.5) -- (-2,-4.5);
            \fill[black] (-2.16,3.35) circle(0.08);
            \fill[black] (-3.79,1.26) circle(0.08);
            \fill[black] (-0.33,4.00) circle(0.08);
            \fill[black] (1.35,3.75) circle(0.08);
            \fill[black] (3.16,2.50) circle(0.08);
            \fill[black] (3.45,-2.00) circle(0.08);
            \fill[black] (1.00,-3.85) circle(0.08);
            \fill[black] (-1.00,-3.85) circle(0.08);
            \fill[black] (-2.33,-3.26) circle(0.08);
            \fill[black] (-3.71,-1.5) circle(0.08);
            \fill[red] (-2,-3.45) circle(0.08);
            \fill[red] (-2,3.45) circle(0.08);
            \fill[red] (2,-3.45) circle(0.08);
            \fill[red] (2,3.45) circle(0.08);
            \draw[arrows=-Latex] (-0.33,4.2) .. controls (-0.8,4.2) and (-1.5,4.05) .. (-2.2,3.6);
            \draw[arrows=-Latex] (-2.3,3.5) .. controls (-2.5,3.6) and (-3.8,2.2) .. (-4,1.4);
            \draw[arrows=-Latex] (1.35,4) .. controls (1.65, 4) and (3, 3.1) .. (3.25,2.7);
            \draw[arrows=-Latex] (-1.00,-4.05) .. controls (-1.3, -4.05) and (-2.25, -3.7) .. (-2.4,-3.4);
            \draw[arrows=-Latex] (-2.5,-3.4) .. controls (-2.8, -3.3) and (-3.6, -2.45) .. (-3.9,-1.75);
            \node at (-2,0) {\Large $S^-$};
            \node at (2,0) {\Large $S^+$};
            \node[red] at (2,3) {\scriptsize $y = 0.3$};
            \node[red] at (-2,3) {\scriptsize $x = 0.9$};
            \node[red] at (-2,-3) {\scriptsize $w = 0.7$};
            \draw[arrows=-Latex] (1.00,-4.05) .. controls (1.5, -4.1) and (3, -3.3) .. (3.6,-2.2);
            \node[red] at (2,-3) {\scriptsize $z = 0.4$};
        \end{tikzpicture}
        \caption{Four local copies of $\mathbb{R}/\mathbb{Z}$ with $x, y, z, w$ and increasing directions labelled.}\label{fig:intervals}
    \end{figure}
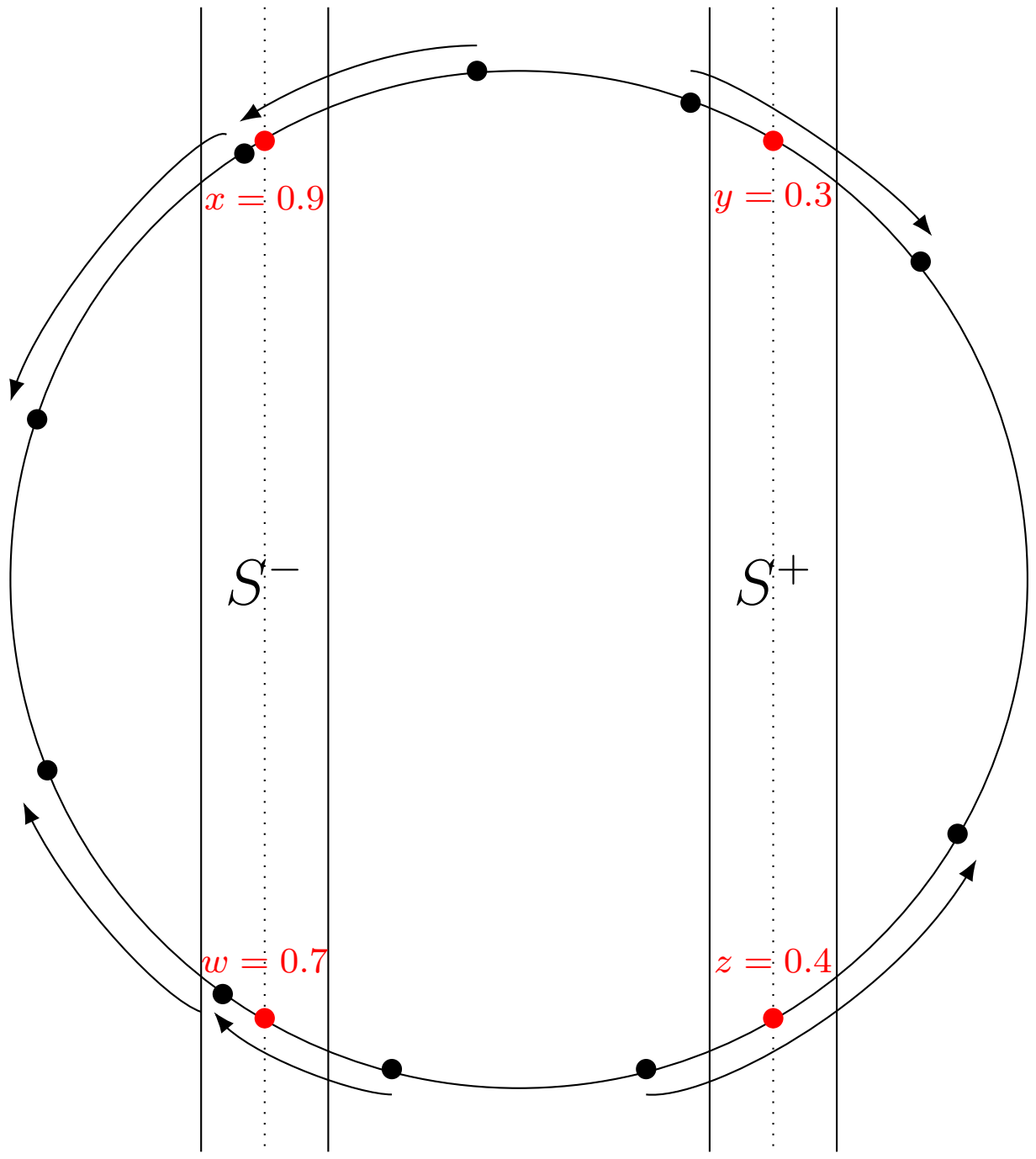
    
    Consider the special gaps as consecutive copies of $\mathbb{R}/\mathbb{Z}$, with the points themselves mapping to $\mathbb{Z}$ and linearly interpolating along the arc in between. Let $I_x$ denote the interval of $\mathbb{R}/\mathbb{Z}$ covered by the side arc, going anticlockwise. If we label its midpoint as $x \in \mathbb{R}/\mathbb{Z}$, we find that
    \[I_x = \left(x - \frac{1}{4\sqrt{3}}, x + \frac{1}{4\sqrt{3}}\right) \pmod{1}.\]
    See Figure \ref{fig:intervals} for an example.

    By construction, we have that $S^- \cap C$ contains a point in the top-left if and only if $0 \in I_x$. Going clockwise, identically define $I_y, I_z, I_w$ as in Figure \ref{fig:intervals} with the directions of increasing position within $\mathbb{R}/\mathbb{Z}$ given as shown. Thus, our goal translates to
    \[0 \not \in I_x \cup I_y \cup I_z \cup I_w.\]
    We have the following interpretations of the two operations in our interval setting. Revolving clockwise uniformly translates
    \[x \uparrow \quad y \downarrow \quad z \uparrow \quad w \downarrow.\]
    Note that $S \cap C$ must stay within the permissible interval. Hence, the maximum translation is the total length of the permissible interval minus the length of the central arc (with $O(\eps^2)$ room for error), namely
    \[L - 2\eps + O(\eps^2) \ge 2\eps + O(\eps).\]
    Converting back to $\mathbb{R}/\mathbb{Z}$, this is a $1/4$ proportion. We also consider the effect of shrinking the radius. Suppose we shrink from $r_1 \mapsto r_2$ with $\Delta := r_2 - r_1$. 

    \begin{figure}[h!]
        \centering
        \begin{tikzpicture}
            \draw [domain=100:140] plot ({6*cos(\x)}, {6*sin(\x)});
            \draw[fill=lightgray] (0,0) -- (-1.5, 2.6) arc [start angle=120, delta angle=3, radius=3cm] -- (0,0);
            \draw [domain=100:140] plot ({5.5*cos(\x)}, {5.5*sin(\x)});
            \fill[black] (0,0) circle(0.08);
            \draw (-3.5,6) -- (-3.5,0);
            \draw (-2.5,6) -- (-2.5,0);
            \draw[red] (0,0) -- (-1.75, 5.73);
            \draw[red] (0,0) -- (-3.6, 4.8);
            \fill[red] (-3, 4.60) circle(0.08);
            \fill[red] (-3, 5.19) circle(0.08);
            \fill[black] (-3.6, 4.8) circle(0.08);
            \fill[black] (-3.3, 4.4) circle(0.08);
            \fill[black] (-1.75, 5.73) circle(0.08);
            \fill[black] (-1.6, 5.26) circle(0.08);
            \draw[dotted] (-3,0) -- (-3,6);
            \draw[dotted] (0,0) -- (-4.5,0);
            \draw[red, dashed] (0,0) -- (-3, 4.60);
            \draw[red, dashed] (0,0) -- (-3, 5.19);
            \node at (-3,1) {\Large $S^-$};
            \node at (-1,2.3) {$\phi$};
            \node at (-0.7,5.4) {$r_2$};
            \node at (-0.71,5.9) {$r_1$};
        \end{tikzpicture}
        \caption{Interpreting the shrinking operation.}\label{fig:shrinkage}
    \end{figure}
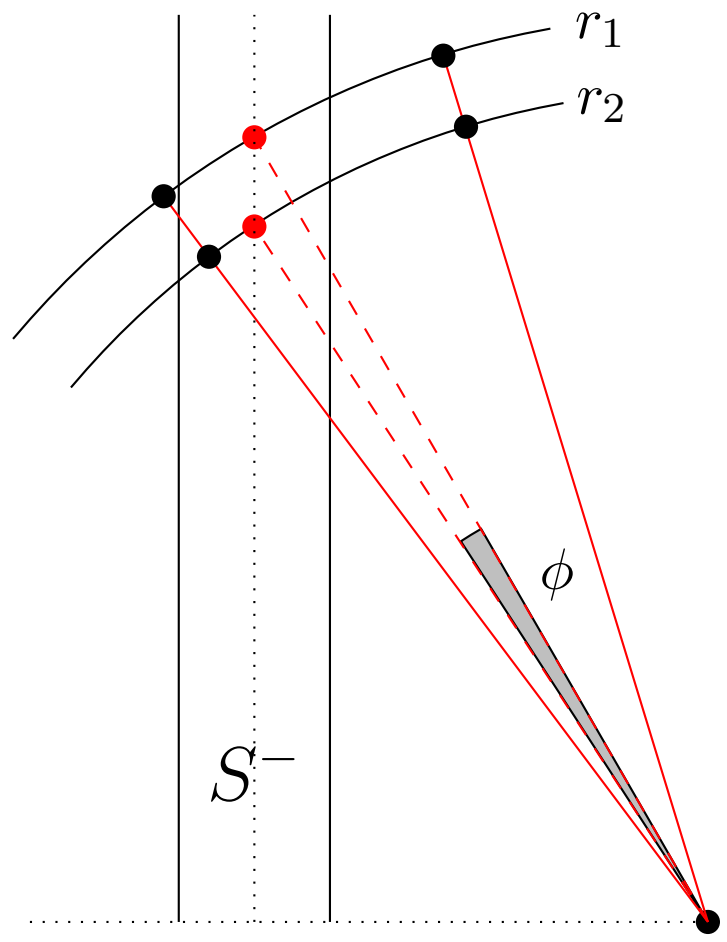

    Let $\phi$ be the angle between the old and new positions of $x$, as shown in Figure \ref{fig:shrinkage}. By the mean value theorem again, there exists $\tilde{r} \in (r_2, r_1)$ such that
    \begin{align*}
        \phi &= \arccos\left(\frac{1}{r_1}\right) - \arccos\left(\frac{1}{r_2}\right) \\
        &= \frac{\Delta}{r_1r_2} \frac{1}{\sqrt{1 - \frac{1}{\tilde{r}^2}}} \\
        &= \left(\frac{1}{2\sqrt{3}} + O(\eps)\right) \Delta.
    \end{align*}
    The arc length lies in $(r_1\phi, r_2\phi)$ hence is $\left(\frac{1}{\sqrt{3}} + O(\eps)\right) \Delta$, since $r_1, r_2 \in (2-\eps, 2-(8\sqrt{3}+1)\eps)$ are both bounded in $\eps$. This quantity only depends on the amount of shrinkage up to a multiplicative $O(\eps)$ error. Hence the effect on the intervals is, up to $O(\eps)$,
    \[x \uparrow \quad y \uparrow \quad z \uparrow \quad w \uparrow.\]
    Note that from the paragraph above, we can safely ignore the central strip. Since we imposed $\Delta \le 8\sqrt{3}\eps$, the increase in each of $x, y, z, w$ is unrestricted. 

    Consider the interval statement: our operations above coincide up to $O(\eps^2)$. Thus, the $\eta$ room for error and lack of restriction in translation let us conclude that for sufficiently small $\eps > 0$, applying our operations yields $0 \not \in I_x \cup I_y \cup I_z \cup I_w$. 
\end{proof}

\subsubsection{Solving the interval problem.} \label{subsec:lower-bound} Without our simplification, we have two sources of error.
\begin{outline}[enumerate]
    \1 The special gaps are long rather than length exactly $8\eps$. In particular, they have length
    \[8\tau \eps + O(\eps^2) \le l \le 8\eps + O(\eps^2).\]
    \1 The permissible interval has length
    \[L \ge 4\tau \eps.\]
\end{outline}
An additional complication is that the two special gaps surrounding (say) the top-left side arc could be significantly different lengths, which results in the endpoints of $I_x$ increasing at different rates with respect to their respective special gaps. Nonetheless, this is encompassed in the following worst-case scenario, which makes edits to the above in turn. 

\begin{outline}[enumerate]
    \1 All special gaps have length $8\tau \eps$. Thus,
    \[|I_x| = |I_y| = |I_z| = |I_w| = \frac{1}{2\sqrt{3}} \tau^{-1}.\]
    Let $\mu := \frac{1}{2\sqrt{3}} \tau^{-1}$.
    \1 Our revolution operation yields a maximum translation of
    \[\frac12 \tau - \frac14.\]
\end{outline}

We claim that we can still solve the relevant interval problem. 

\begin{prop}\label{prop:interval}
    Suppose $0 < \delta < \delta_0$, where
    \[\delta_0 = \frac{1}{259}.\]
    Let $\tau := \frac{1-25\delta}{1+\delta}$ and $\mu := \frac{1}{2\sqrt{3}}\tau^{-1}$. 
    
    Let $x, y, z, w \in \mathbb{R}/\mathbb{Z}$. Let $I_x = (x - \frac{\mu}{2}, x + \frac{\mu}{2}) \pmod{1}$ and similarly define $I_y, I_z, I_w$. Then one can apply the following operations until $[0, 0.01] \cap (I_x \cup I_y \cup I_y \cup I_z) = \varnothing$.
    \begin{outline}[enumerate]
        \1 Increase $x, z$ and decrease $x, w$ uniformly by $t \in [0, \frac12 \tau - \frac14]$.
        \1 Increase $x, y, z, w$ uniformly by $t \in [0,1)$.
    \end{outline}
\end{prop}
\begin{proof}
    We can remove the second operation to change the goal to excluding some sub-interval of length $\eta$. Indeed, then we could translate so the left endpoint is 0. 

    Merge $A = I_x \cup I_z$ and $B = I_y \cup I_w$ since they move at the same velocities within a pair, and opposite to the other pair. In this step, we may without loss of generality assume $I_x \cap I_z = \varnothing$ and $I_y \cap I_w = \varnothing$. In particular, we let $A, B$ be disjoint unions of open intervals of length $\mu$. 

    Furthermore, we are allowed to fix a frame of reference where
    \[A = (0, \mu) \cup (a, a + \mu), \quad \mu \le a \le 1-\mu,\]
    \[B = (b_1, b_1 + \mu) \cup (b_2, b_2 + \mu), \quad 0 \le b_1 < b_1 + \mu \le b_2.\]
    Since $A$ is stationary, $B$ moves at most $\tau - \frac12$. 

    Consider setting $\tau > 0.9$, then suppose for the sake of contradiction that $A \cup B$ always meets every interval of length 0.01. 

    \under{Case 1}: We have $a \le \mu+0.01$ and $b_2 \le b_1 + \mu + 0.01$. Without loss of generality, $B$ moves leftwards: if not, we can swap the roles of $A, B$. 
    
    To initially meet $[a_1+\mu, a_1+\mu+0.01]$, we require $b_1 \le a+\mu + 0.01$. To ensure $[0.99, 1]$ is always intersected, we also need the right endpoint of $B$ to always lie past $0.99$. In particular, after translating the maximal amount, we find that
    \[(b_2 + \mu) - \left(\tau - \frac12\right) \ge 0.99.\]
    Combining the inequalities together, we obtain $b_2 \le 3\mu + 0.03$. Thus, we have
    \[\frac{2}{\sqrt{3}}\tau^{-1} - \tau \ge 0.46.\]
    This gives $\tau \le 0.87$, a contradiction.

    \under{Case 2}: Without loss of generality, suppose $a > \mu + 0.01$. Again, suppose $B$ moves leftwards: if not, shift $A$ to the other endpoint of $\mathbb{R}/\mathbb{Z}$ and repeat the argument below but reflected about 1/2. 

    Since $\tau > 0.87$, we have $\tau - \frac12 > \mu$. Thus, we can translate $B$ until $b_1 = 0$. At this instant, we would have
    \[A \cup B = (0, \mu) \cup (a, a+\mu) \cup (b_2, b_2 + \mu).\]
    The complement within $\mathbb{R}/\mathbb{Z}$ is at most three (closed) intervals of total length at least $1 - 3\mu$ by the union bound. Since none of the three intervals can have length more than 0.01, we obtain $1 - 3\mu \le 0.03$. This rearranges to $\tau \le 0.9$, giving a contradiction. \\

    Thus, if
    \[\frac{1-25\delta}{1 + \delta} > 0.9\]
    then the conclusion of the proposition holds. In particular, we obtain
    \[\delta < \delta_0 := \frac{0.1}{25.9} \approx 0.00386.\qedhere\]
\end{proof}

We can now prove our lower bound improvement.

\begin{proof}[Proof of Theorem \ref{thm:lower}]
    Take the alternative perspective in Section \ref{subsec:perspective}. Thus, we may assume that consecutive $\alpha_j$ forming $\Theta$ have angles differing by at most $2\pi/x_0$. Then we apply Lemma \ref{lmm:avg-looking} to identify a revolution of our points such that all special gaps are long. 

    Combining the worst-case edits made at the start of this section with the simplification given by Lemma \ref{lmm:reduction}, we find that if Proposition \ref{prop:interval} holds, then there exists a choice of revolution and radius such that none of our points lie in $E(\eps)$. Since Proposition \ref{prop:interval} does hold, we find $z \in \mathbb{C}$ not covered by the rotations of the pyjama stripe by $\Theta$. 
\end{proof}

\section{Concluding remarks}

Our upper bound in Theorem \ref{thm:pyjamas} and lower bound in Theorem \ref{thm:lower} still differ by a quasipolynomial gap, notably with the lower bound still linear in $\eps^{-1}$. Even with stronger hypothetical results on $M_k$ covering $(\mathbb{Z}[i]/Q\mathbb{Z}[i])^\times$ for $N(Q) \le \exp(Ck)$, where $C > 0$ is an absolute constant, the naive generalisation of our proof still yields quasipolynomial bounds, as an artifact of the growing space. Indeed, in the language of Section \ref{sec:minor}, we require $\eta$ to decay (at least) exponentially with $n$, and $k$ to grow at least super-linear in $n$, which yields an $\exp(n^2)$ upper bound for $N_0$. In line with Green's proposal \cite{green2024}, a polynomial upper bound seems believable. 

One can consider generalising our lower bound approach by incorporating circles (or balls) of larger radius. It is plausible that analysing an infinite sequence of circles and their neighbourhoods provides a noticeably stronger result. We chose circles to exploit spherical symmetry; it is also possible to experiment with other compact subsets of the plane. \\

We also pose the following generalisation of the Pyjama Problem, phrased to resemble the original Problem \ref{prob:pyjama}. 

\begin{problem}\label{prob:shrinking}
    Let $(\eps_n)_{n = 1}^\infty$ be a sequence of positive real numbers with divergent sum. Does there exist a finite set $\{\psi_1, \dots, \psi_N\}$ of unit-norm complex numbers such that
    \[\bigcup_{n = 1}^N E_{\psi_n}(\eps_n) = \mathbb{C}?\]
\end{problem}

The Pyjama Problem asks this for the constant sequences $(\eps)_{n = 1}^\infty$, where $\eps > 0$. 

Certainly an infinite number of rotations suffices. Indeed, we can consider rotations of the central strip only. Given any $R > 0$ and $m \in \mathbb{N}$, since $\sum_{n = m}^\infty \eps_n$ diverges, there exists some $k \in \mathbb{N}$ such that
\[\sum_{n = m}^{m+k} \eps_n \ge 2\pi R,\]
and in particular, we can greedily cover $S_R(0)$. This also covers the entirety of $\overline{B}_R(0)$. Taking the union of these rotations for $R = 1, 2, \dots$ then covers the plane. 

A particularly interesting instance of Problem \ref{prob:shrinking} is the harmonic sequence, given by $\eps_n = 1/n$.

\subsubsection*{AI statement.} At the start of this project to improve bounds on the Pyjama Problem, James Leng communicated to the author the progress that ChatGPT 5.6 Pro could make. For the lower bound, it produced no ideas. For the upper bound, it could not make any progress using a growing set of primes, but proposed the reduction to continuous characters to shave an exponential.

We have re-written the output ourselves to ensure clarity, and contained it entirely within Section \ref{sec:pigeon}. All the ideas in the other sections were discovered and written-up by the author without AI assistance. 

\section{Acknowledgments}

The author is incredibly grateful to Noah Kravitz and James Leng for introducing the problem and the countless helpful discussions since. 

This research started as part of the Summer Research in Mathematics (SRIM) Programme organised by the Department of Pure Mathematics and Mathematical Sciences at the University of Cambridge. During this period, the author was supervised by Timothy Gowers and is very grateful for the meet-ups and mentorship. The author would also like to thank Joni Teräväinen, Cédric Pilatte and Péter Varjú for many insightful conversations. This was made possible by the Trinity College Summer Studentship Scheme, which supported the author throughout this period. 

Part of this project was also carried out while the author was visiting Emory University. The author is grateful to Cosmin Pohoata for their hospitality and thanks the Global Talent Fund for supporting the stay in Atlanta. 

\bibliographystyle{plain}
\bibliography{citations}

\end{document}